\documentclass[psamsfonts]{amsart}

\usepackage{latexsym}
\usepackage{amscd}
\usepackage{amssymb}
\usepackage{amsfonts}
\usepackage{amsmath}
\usepackage{xy}
\xyoption{all}
\usepackage{amsthm}
\usepackage{enumerate}
\usepackage{hyperref}
\usepackage{eprint}
\usepackage{epsfig}
\usepackage{graphics}

\newtheorem{thm}{Theorem}[section]
\newtheorem{lem}[thm]{Lemma}
\newtheorem{defn}[thm]{Definition}
\newtheorem{prop}[thm]{Proposition}
\newtheorem{quoteprop}[thm]{``Proposition''}

\newtheorem{eg}[thm]{Example}

\newtheorem{rmk}[thm]{Remark}

\numberwithin{equation}{section}

\def\a{\alpha}
\def\b{\beta}

\def\R{\mathbb{R}}
\def\H{\mathbb{H}}

\def\C{\mathbb{C}}
\def\Q{\mathbb{Q}}
\def\Z{\mathbb{Z}}

\def\P{\mathbb{P}}
\def\N{\mathbb{N}}
\def\del{\partial}

\def\ra{\rightarrow}

\def\del{\partial}

\def\Ainf{A_{\infty}}
\def\FS{\mathcal{FS}}
\def\Fuk{\operatorname{Fuk}}
\def\Coh{\operatorname{Coh}}
\def\Qcoh{\operatorname{QCoh}}
\def\Perf{\operatorname{Perf}}
\def\hol{\operatorname{Hol}}

\def\Pic{\operatorname{Pic}}
\def\Int{\operatorname{Int}}

\def\mod{\operatorname{mod}}
\def\Mod{\operatorname{Mod}}
\def\deg{\operatorname{deg}}

\def\id{\operatorname{id}}
\def\Id{\operatorname{Id}}
\def\RHom{\operatorname{RHom}}
\def\Hom{\operatorname{Hom}}

\def\Ext{\operatorname{Ext}}
\def\Int{\operatorname{Int}}
\def\Out{\operatorname{Out}}
\def\Hol{\operatorname{Hol}}

\begin{document}

\author{Matthew Robert Ballard}
\address{Department of Mathematics, University of Washington,
Seattle, WA 98195, USA}
\email{ballard@math.washington.edu}


\title{Meet Homological Mirror Symmetry}
\maketitle

\begin{abstract}
 In this paper, we introduce the interested reader to homological mirror symmetry. 
After recalling a little background knowledge, we tackle the simplest cases of 
homological mirror symmetry: curves of genus zero and one. We close by outlining 
the current state of the field and mentioning what homological mirror symmetry 
has to say about other aspects of mirror symmetry.
\end{abstract}

\tableofcontents

\section{Introduction}

As it stands today, mirror symmetry is a well-established field of mathematics. Seventeen years have passed since physicists correctly predicted the number of rational curves of given 
degrees on the Fermat quintic \cite{COGP91}. In the intervening years, the physical inspiration 
has energized the mathematics community. In many cases, including anti-canonical hypersurfaces and Calabi-Yau complete intersections in toric varieties, mathematical mirror constructions have been realised allowing one to see the mirror reflection in the Hodge diamond and establish the relation between curve counts and period integrals. Unfortunately, rigorous formulations of the 
physical arguments yielding these predictions have not arisen. Such rigorous formulations would allow mathematicians to tap into the deep physical understanding of the situation, but modern physical tools, such as the path integral, lie just beyond the current scope of mathematics. Consequently, to fully understand the phenomenon of mirror symmetry, mathematicians must follow the physical inspiration and develop their own overarching approaches. 

At the International Congress of Mathematicians in 1994, M. Kontsevich proposed such an approach 
\cite{Kon95}. As is it commonly known now, homological mirror symmetry, or the homological 
mirror conjecture, reformulates mirror symmetry as an equivalence of triangulated categories built 
from different aspects of the Calabi-Yau geometry of two manifolds, called mirror manifolds. The Calabi-Yau metric 
is determined by two pieces of data --- the complex structure and the K\"ahler form. Each piece 
lands in a different branch of mathematics. If we only care about the 
complex structure, we are naturally led into algebraic geometry. If we instead remember the 
K\"ahler form, we study symplectic geometry. Upon passage to the mirror manifold, we should see 
an exchange of these two geometries. Kontsevich sought to capture this exchange. Today, homological mirror symmetry remains an intriguing and daunting challenge to mathematics. 
The scope has expanded to encompass other manifolds beyond Calabi-Yau manifolds. It has become a 
powerful source of inspiration motivating insights in algebraic geometry, symplectic geometry, 
homological algebra, noncommutative geometry, and beyond.

Mixing all these fields, homological mirror symmetry is a very attractive conjecture, but it 
remains outside the common knowledge of the working mathematician. Why? Along with its relative youth, its proper formulation requires an imposing amount of unfamiliar technology. This article offers a  introduction to homological mirror symmetry through two explicit examples. The cases of $\P^1$ and elliptic curves are very concrete. We give a holistic approach that treats both sides of mirror symmetry for $\P^1$, instead of speaking to one side without reference to the other as is often the case in the literature. Since homological mirror symmetry manifests itself in noticeably different ways whether one considers Fano or Calabi-Yau varieties, we also review the case of elliptic curves. Seeing both cases side by side will hopefully give the reader a deeper appreciation of the duality.
To move into examples we need to review a little of the necessary formalism. This is done in section 
\ref{sec:commonknowledge}. With this over, we move on to the examples in sections \ref{sec:genuszero} 
and \ref{sec:genusone}. Algebraic curves of genus zero and one provide simple case studies. Within 
these, the reader can meet concrete incarnations of the relevant categories and appreciate the unexpected 
equivalences that homological mirror symmetry predicts. After the examples have been covered, we 
outline the current state of knowledge in the field in section \ref{sec:furtherresults} and then mention 
how homological mirror symmetry relates to other aspects of mirror symmetry in section \ref{sec:HMSMS}.

The author heartily thanks his adviser, Charles Doran, for encouraging the creation of this survey, the referee for his careful reading, and Ursula Whitcher for her time and numerous suggestions to improve the article. All errors belong solely to the author.

\section{Building some common knowledge}
\label{sec:commonknowledge}

Part of the difficulty in dealing with homological mirror symmetry is the breadth of knowledge 
required for a proper formulation. Before we can dive into the promised simple examples, we recall 
some terminology and results from homological algebra, algebraic geometry, and symplectic geometry.

\subsection{Homological algebra}

I assume the reader has some familiarity with derived categories, at least in the case of modules over an associative algebra, Hochschild cohomology, and dg-algebras. A good reference for homological algebra is \cite{GM03}, a good reference for homological algebra and Hochschild cohomology is \cite{Wei94}, and a good reference for dg-algebras is \cite{Kel06}.

Our goal is here is to define a triangulated category which appears on each side of mirror symmetry. The main algebraic tool is the $\Ainf$-algebra. $\Ainf$-algebras are avoidable in algebraic geometry but not in symplectic geometry. (However, the $\Ainf$-algebras that appear in homological mirror symmetry for $\P^1$ are honest associative algebras).

\begin{defn}
 An $\Ainf$-algebra over a base field $k$ is a graded $k$-module $A$ with
$k$-linear maps $m_n: A^{\otimes n} \ra A$ of degree $2-n$ for each $n > 0$,
satisfying the following quadratic relations for all $n>0$:
\begin{displaymath}
 \sum_{r,s} (-1)^{rs+(n-r-s)} m_l(\id^{\otimes r} \otimes m_s \otimes \id^{\otimes n-r-s}) = 0
\end{displaymath}
\end{defn}
Let us look at the first three of these relations. The first relation says $m_1^2=0$ so $(A,m_1)$ is a chain complex. The second 
says $m_2$ is a chain map when we use the differential $1 \otimes m_1 + m_1 \otimes 1$ on $A^{\otimes 2}$. 
The third says that $m_2$ is associative up to a homotopy $m_3$. Thus, we can pass to the cohomology 
with respect to $m_1$ and get an associative algebra $H(A)$.

Familiar examples of $\Ainf$-algebras are associative algebras, where $m_n=0$
for $n\not =2$, and dg-algebras, where $m_n = 0$ for $n\not = 1,2$. An $\Ainf$-algebra with $m_1=0$ is called minimal.

We can (and often have to) do something troublesome and add a degree $0$ multiplication $m_0: k \ra A$
and continue to require the quadratic relations to hold.
For example, the first two relations become
\begin{displaymath}
 \begin{aligned}
  & m_1(m_0) = 0
  & m_1(m_1) + m_2(m_0,\id) + m_2(\id,m_0) = 0
 \end{aligned}
\end{displaymath}
This in general destroys the possibility of taking cohomology. An $\Ainf$-algebra possessing such 
an $m_0$ is called curved or obstructed. If we assume, that $m_0$ is
central (with respect to $m_2$) and $m_n(\id^{\otimes r} \otimes m_0 \otimes
\id^{\otimes n-r-1}) = 0$ for all $n > 2$ we can once again take cohomology 
$H(A)$. In this case, $A$ is called weakly obstructed. 
Below our $\Ainf$-algebras will be unobstructed unless explicitly indicated.

\begin{defn}
 A morphism $f: A \ra B$ of $\Ainf$-algebras over $k$ is a collection of $k$-linear maps 
$f_n: A^{\otimes n} \ra B$ of degree $1-n$ satisfying
\begin{displaymath}
 \sum (-1)^{\heartsuit} m_k(f_{i_1}\otimes \cdots \otimes f_{i_k}) = \sum (-1)^{sr+(n-s-r)} f_l(\id^{\otimes s} \otimes m_r \otimes \id^{\otimes n-r-s})
 \label{ainfmorphismrelations}
\end{displaymath}
where $\heartsuit = (k-1)(i_1-1) + (k-2)(i_2-1) + \cdots + (i_{k-1}-1)$.
\end{defn}
The first relation says that $f_1$ commutes with the differential. The second says
$f_1$ respects $m_2$ up to $f_2$. $f$ is a called a quasi-isomorphism
if $f_1:H(A) \ra H(B)$ is an isomorphism.

\begin{rmk}
 These definitions become more compact and perhaps clearer when one passes from $A$
to the bar complex $B(A)$ on $A$. The $\Ainf$-algebra structure on $A$ is equivalent to a 
coderivation on $B(A)$. $\Ainf$-algebra morphisms are then coalgebra morphisms commuting 
with the coderivations. For a reference, see \cite{GJ90}.
\end{rmk}

Now we recall a result that allows us to pass to cohomology of an $\Ainf$-algebra without losing 
information.

\begin{thm}
 \cite{Kad80,Mer99,Mar06} Given an $\Ainf$-algebra $A$ over a field $k$, choose a splitting $A = H \oplus B \oplus D$ over $k$, where $H$ is the cohomology and $m_1: D \ra B$ is an isomorphism. 
Then, there exists an $\Ainf$-algebra structure on $H(A)$ with zero first order composition, second order composition induced by $m_2$, and
quasi-isomorphisms $i: H(A) \ra A, \pi: A \ra H(A)$. Moreover, $\pi_1$ is the associated restriction 
map $A \ra H$ and $i_1$ is the associated inclusion map $H \ra A$.
\label{thm:minimalmodel}
\end{thm}

Consider for a moment a dg-algebra $A$. Applying this procedure yields a minimal $\Ainf$-algebra 
structure on $H(A)$. Thus, the chain-level data lost when we just consider $H(A)$ has been 
transmogrified, returning as the higher compositions in the $\Ainf$-structure.

How complex can $\Ainf$-algebras be? A partial answer to this question is given by the following observation.

From theorem \ref{thm:minimalmodel}, there is an $\Ainf$-structure on $H(A)$ which 
makes it quasi-isomorphic to the $\Ainf$-algebra $A$. Given a minimal 
$\Ainf$-algebra $(A,m_n)$, take $m_k$ to be the first non-zero operation with $k>2$. Then, the 
first $\Ainf$-relation involving $m_k$ is
\begin{displaymath}
\sum_r (-1)^r m_k(\id^{\otimes r} \otimes m_2 \otimes \id^{\otimes k-1-r}) + m_2(m_k,\id) - m_2(\id,m_k) = 0
\end{displaymath}
This equation states that $m_k$ is a Hochschild cochain for 
the algebra $(A,m_2)$. Suppose we wanted to find another minimal $\Ainf$-structure $m_n'$ on $(A,m_2)$ 
which is isomorphic to our original one. Then, we have our collection $f_k: A^{\otimes k} \ra A[1-k]$. We need $f_1$ to be 
an automorphism of the algebra $(A,m_2)$, consequently we can inductively solve
\begin{displaymath}
 \sum (-1)^{\heartsuit} m_k(f_{i_1}\otimes \cdots \otimes f_{i_k}) = \sum (-1)^{s+r(n-s-r)} f_l(\id^{\otimes s} \otimes m_r \otimes \id^{\otimes n-r-1})
\end{displaymath}
for $m_n'$. Thus, the $f_n$ and $m_n$ uniquely determine $m_n'$. So, if we want to find an isomorphic 
$(A,m_n')$ with $m_l' = 0$ for $2<l\leq k$, we just need to solve
\begin{displaymath}
 \sum_s f_l(\id^{\otimes s} \otimes m_2 \id^{\otimes n-r-1}) + m_2(f_l,\id) + m_2(\id,f_l) = m_k(f_1^{\otimes k})
\end{displaymath}
and then take $f_l=0$ for $1<l<k$ and $f_l$ arbitrary for $l>k$. Applying an automorphism, we can set 
$f_1 = \id$. The previous equation says that $m_k$ is a 
Hochschild coboundary for $(A,m_2)$. Thus, we get the following result. 

\begin{lem}
 If the relevant pieces of the Hochschild cohomology of 
$(A,m_2)$ are zero, we can trivialise up to any arbitrary order. From the construction, the limiting composition of these maps exists, giving a trivialisation for the $\Ainf$-structure.
\label{lem:HochtrivAinf}
\end{lem}

As in the case of an associative algebra, we wish to form a category of modules over an $\Ainf$-algebra $A$. We do so in a manner strongly analogous to forming chain complexes of modules over an associative algebra.

\begin{defn}
 A (right) module $M$ over an $\Ainf$-algebra $A$ is a graded module over $k$ equipped with $k$-linear maps $m^M_n: A^{\otimes n-1} \otimes M \ra M$ of degree $2-n$ for each $n > 0$: satisfying the following quadratic relations for all $n>0$:
\begin{gather*}
 \sum_{r,s,n-r-s>0} (-1)^{rs+(n-r-s)} m^M_l(\id^{\otimes r} \otimes m_s \otimes \id^{\otimes n-r-s}) + \\ \sum_{u,v} (-1)^{uv} m^M_l(\id^{\otimes u} \otimes m_v^M)   = 0
\end{gather*}
 A morphism $g: M \ra N$ of $A$-modules is a collection of $k$-linear maps $g_n: A^{\otimes n-1} \otimes M \ra N$ satisfying quadratic relations similar to the case of morphisms of $\Ainf$-algebras. For more details, see \cite{GJ90}.
\end{defn}

Let us repackage this definition. Instead of thinking of $A$ as an $\Ainf$-algebra, we shall think of it as an $\Ainf$-category with one object $*$. The morphisms in this category are only the endomorphisms of $*$ and $\Hom(*,*) := A$. To give a category $\mathcal{A}$ an $\Ainf$-structure, we need to define multi-compositions of morphisms 
\begin{displaymath}
 m_n: \Hom(X_0,X_1) \otimes \cdots \otimes \Hom(X_{n-1},X_n) \ra \Hom(X_0,X_n)
\end{displaymath}
satisfying the $\Ainf$-relations. Here $X_i$ are objects of $\mathcal{A}$. In our case, we just use the operations $m_n: A^{\otimes n} \ra A$ coming from the $\Ainf$-algebra structure of $A$. 

For an $\Ainf$-category $\mathcal{A}$, each object $X$ of $\mathcal{A}$ furnishes a functor
$\Hom_{\mathcal{A}}(X,\cdot)$ from $\mathcal{A}$ to chain complexes over $k$, $Ch(k)$. These are naturally $\Ainf$-functors (for more, see \cite{Kel01}).

\begin{defn}
 Given two $\Ainf$-categories $\mathcal{A}$ and $\mathcal{B}$. An $\Ainf$-functor $\mathcal{F}: \mathcal{A} \ra \mathcal{B}$ is an assignment of objects $X \mapsto \mathcal{F}(X)$ and a collection of maps
\begin{displaymath}
 \mathcal{F}_n: \Hom_{\mathcal{A}}(X_0,X_1) \otimes \cdots \otimes \Hom_{\mathcal{B}}(X_{n-1},X_n) \ra \Hom_{\mathcal{B}}(\mathcal{F}(X_0),\mathcal{F}(X_n))
\end{displaymath}
which satisfy quadratic relations similar to those given in the definition of a morphism of $\Ainf$-algebras.
\end{defn}

The 
$\Ainf$-version of the Yoneda embedding is given by sending $X$ to $\Hom_{\mathcal{A}}(X,-)$. This gives an $\Ainf$-functor $Y$ from  $\mathcal{A}$ to the category of $\Ainf$-functors from $\mathcal{A}^{op}$ to $Ch(k)$. Denote this category by $\Mod \mathcal{A}$. The standard Yoneda embedding is full and faithful. After taking $H^0$, 
$Y$ reduces to the standard Yoneda embedding. Since $H^0(Y)$ is an equivalence onto its image, we say that $Y$ is a quasi-equivalence onto its image.

If we consider the $\Ainf$-category with a single object and morphisms algebra $A$, this gives another definition of $A$-module. It is a good exercise for the reader to translate between the two definitions of $A$-modules for an $\Ainf$-algebra. 

Since any $A$-module is, in particular, a chain complex, we have a notion of quasi-isomorphism in the category of $A$-modules. In analogy with the case of an associative algebra, we wish to invert these. When we do, we get $D(A)$, the derived category of $A$-modules. As in the case of associative algebras, the resulting category is triangulated. We let $D^{\pi}(A)$ be the smallest triangulated subcategory of $D(A)$ containing $A$ and closed under taking triangles, direct sums, and direct summands. To finish, we recall the following useful result. For a proof, see \cite{SeiDR}, where the reader can find more details about the construction of $D(A)$.

\begin{prop}
 If $A$ and $B$ are quasi-isomorphic $\Ainf$-algebras, then $D(A)$ is triangle equivalent to $D(B)$ and 
$D^{\pi}(A)$ is triangle equivalent to $D^{\pi}(B)$.
\label{prop:quasiequivcat}
\end{prop}

\subsection{Algebraic geometry}

Algebraic geometry is the most natural geometric field for the application of homological algebra. 
Indeed, homological algebra permeated algebraic geometry long ago. We shall only recall a small amount. For more details, see \cite{Har66}.

Let us recall some standard 
abelian categories associated to an algebraic variety $X$. Any algebraic variety $X$ comes with 
its sheaf of functions $\mathcal{O}_X$. A sheaf $\mathcal{E}$ on $X$ is called a quasi-coherent sheaf 
if there is an action of $\mathcal{O}_X$ on $\mathcal{E}$ and, locally in the Zariski topology, 
$\mathcal{E}$ is the cokernel of a morphism between free $\mathcal{O}_X$-modules. $\mathcal{E}$ is 
called coherent if locally, in the Zariski topology, $\mathcal{E}$ is the cokernel of a morphism 
between free, finite-rank $\mathcal{O}_X$-modules. Let us restrict ourselves to the case that $X$ is quasi-projective 
over a field $k$. In this case, the category of coherent sheaves on $X$, $\Coh(X)$, 
is an abelian category. The category of quasi-coherent sheaves, $\Qcoh(X)$, is also abelian and it 
possesses enough injective objects. Therefore, we can resolve any quasi-coherent sheaf by a bounded below complex of injective sheaves. Thus, we can form the derived category of $\Coh(X)$, or $\Qcoh(X)$, by taking the homotopy category of the subcategory of all bounded below complexes of injectives with bounded coherent cohomology, or by taking the homotopy category of the subcategory of all injectives with quasi-coherent cohomology.

An important class of coherent sheaves is the locally free coherent sheaves, i.e. ones which are locally 
isomorphic to a finite-rank free $\mathcal{O}_X$-module. Given a locally free coherent sheaf $\mathcal{E}$,
we can associate to it an algebraic vector bundle $E$ whose sheaf of sections is the dual sheaf 
of $\mathcal{E}$, i.e. the sheaf whose sections over $U$ are given by 
$\Hom_U(\mathcal{E}(U),\mathcal{O}_X(U))$. This gives a contravariant equivalence between the 
categories of locally free coherent sheaves and algebraic vector bundles.

\subsection{Symplectic geometry background}

In this section, we review some of the underlying notions of symplectic geometry. With this knowlege in hand, we proceed, in this section, to define Lagrangian intersection Floer homology. An excellent reference for all things symplecto-topological is \cite{MS98} and, 
correspondingly, an excellent introduction to the uses of $J$-holomorphic curves in symplectic topology is \cite{MS04}

Let $M$ be a smooth manifold and $\omega$ an anti-symmetric two-form.

\begin{defn}
$\omega$ is a symplectic form if $d\omega = 0$ and $\omega$ is
non-degenerate, i.e. the pairing on vectors in $T_xM$ at all points $x \in M$ is
non-degenerate.
\end{defn}

Note that a symplectic manifold must have even dimension.

\begin{eg}

\begin{enumerate}
\item The canonical example of a symplectic manifold is the cotangent space
$T^*X$ to any smooth manifold $X$. Let us look at an open chart $U$ and denote
the coordinates on $U$ by $p$ and denote the coordinates in the fiber direction
by $q$. $T^*U \cong U \times \R^n$. We let $\sigma_{can} = q \ dp$. One can
check that even though we have defined this locally, it is globally
well-defined. We set $\omega_{can} = d\sigma_{can} = dq \wedge dp$. $\omega_{can}$ is
clearly non-degenerate as vectors in the base pair with vectors in the fiber.
$(T^*X,\omega_{can})$ is also an important example of an exact symplectic manifold, that is, a symplectic manifold for which $\omega = d \sigma$ for some one-form $\sigma$.
\item Any smooth projective complex algebraic variety is a symplectic manifold.
The Fubini-Study K\"ahler form on $\P_{\C}^n$ restricts to a symplectic form on 
any complex submanifold.
\end{enumerate}
\end{eg}

The following result shows that symplectic manifolds are all locally
isomorphic.

\begin{prop}
 (Darboux) Given a symplectic manifold $(M^{2n},\omega)$ and a point $x \in M$, then there is 
a neighborhood $U$ of $x$ and a set of coordinates $$(x^1,\ldots,x^n,y^1,\ldots,y^n)$$ on $U$ such 
that $\omega|_U = \sum_i dx^i \wedge dy^i$.
\end{prop}

If $\omega$ is an anti-symmetric two-form on a vector space $T$, for any
subspace $S$ we can define the symplectic orthogonal $S^{\perp} = \{v \in T :
\omega(s,v) = 0 \ \forall s \in S\}$. If $\omega$ is non-degenerate, then $S
\oplus S^{\perp} = T$. A subspace $S$ is isotropic if $S \subset
S^{\perp}$ and coisotropic if $S^{\perp} \subset S$. It is called Lagrangian if
it is both isotropic and coisotropic and symplectic if $S \cap S^{\perp} =
\{0\}$. Given an symplectic manifold $(M,\omega)$, we say that an embedded
submanifold $S$ is isotropic (resp. coisotropic, Lagrangian, symplectic) if all
the tangent spaces are isotropic (resp. coisotropic, Lagrangian, symplectic).

For $T^*X$ with $\omega_{can}$, the zero section $X$ is a Lagrangian submanifold,
as is each fiber. The general case is similar.

\begin{prop}
(Weinstein) Let $L$ be a Lagrangian submanifold of a symplectic manifold
$(M,\omega)$. Then there is a tubular neighborhood $U$ of $L$ which is
symplectomorphic to a neighborhood of $L$ in $T^*L$ with $\omega_{can}$.
\end{prop}

Since a symplectic form is non-degenerate, it affords a way to turn vectors into
one-forms and vice versa. If $\eta$ is a one-form, we can define $\eta^{\vee}$
to be the vector field such that $\omega(\eta^{\vee},-) = \eta$. We often start with the differential of a smooth function $f: M \ra \R$,
or of a time-varying function $f: M \times \R \ra \R$. Given such a function, we
define the Hamiltonian vector field $X_{f_t}$ to be $df_t^{\vee}$. Given a
vector field, we can often integrate it out to a diffeomorphism called a flow
$\phi^t: M \ra M$. The flow associated to $X_{f_t}$ is called a Hamiltonian
flow. One can check that since $df_t$ is closed, Hamiltonian flows preserve the
symplectic form. Note that given a submanifold $N$ of $M$, $\phi^t(N)$ defines
an isotopy between $N$ and $\phi^1(N)$. This is called a Hamiltonian isotopy.

\begin{defn}
A bundle endomorphism $J: TX \ra TX$ is called an almost complex structure if $J^2 =
-1$.
\end{defn}

The simplest example of an almost complex structure is multiplication by
$\imath$ on $\C \cong \R^2$. An almost complex structure on a vector space is
the same thing as a complex structure. Complex manifolds therefore have
natural almost complex structures, namely multiplication by $\imath$. An
almost complex structure $J$ that makes $X$ into a complex manifold is called a
complex structure. However, not all almost complex structures are complex
structures. The $\imath$-eigenspaces of $J$ on $TX \otimes \C$ must be
involutive. See \cite{NN57}.

\begin{defn}
If $(M,\omega)$ is a symplectic manifold, then an almost complex structure $J$
is called compatible if $\omega(J-,-)$ is a Riemannian metric, i.e. symmetric,
positive definite, and non-degenerate.
\end{defn}

Given a symplectic manifold $(M,\omega)$, we can always
locally find a compatible almost complex structure using Darboux's theorem. Since Riemannian metrics
form a convex space (we can add them), we can find a partition of unity and
patch together our local compatible almost complex structures into a global
almost complex structure $J$. Given a complex manifold $Y$ with complex
structure $j$, we say that a map $\phi: Y \ra M$ is $J$-holomorphic if $d\phi
\circ j = J \circ d \phi$. If $Y = M = \C$ and $j = J = \imath$, then this
equation is equivalent to the set of Cauchy-Riemann equations. Solutions of the
Cauchy-Riemann equations satisfy nice properties (all of complex algebraic
geometry stems from this.) Gromov realised \cite{Gro85} that $J$-holomorphic
maps have many of the same nice properties. He introduced the study of $J$-holomorphic 
curves into symplectic manifolds and gave a useful new tool to symplectic geometry.

\subsection{Floer cohomology and Fukaya categories}

Now we can outline Lagrangian intersection Floer cohomology. Let $L$ and $L'$ be
two closed Lagrangian submanifolds of a compact  symplectic manifold. If they
are not transverse, replace $L$ by a Hamiltonian isotopic Lagrangian. We can
then assume that $L$ and $L'$ intersect transversely and thus $L \cap L'$ is a
finite set. Given two intersection points $p$ and $q$, consider the set of
$J$-holomorphic maps $\phi: D \ra M$ from the unit disc with two marked points
at $-1$ and $1$ such that $\phi(\del D \cap \H) \subset L, \phi(\del D \cap
(\C-\H)) \subset L', \phi(1) = q,$ and $\phi(-1) = p$. The expected dimension of the space of solutions can be determined as follows.
The pullback tangent bundle $\phi^*TM$ is trivial since we are working with a
disc. There is a real sub-bundle along the boundary determined by the tangent
spaces to the Lagrangians. We can change our trivialisation so that
$\phi^*(T_pL') = \imath \phi^*(T_pL)$ and $\phi^*(T_qL) = \imath \phi^*(T_qL')$.
If we rotate at $p$ and $q$, we get a loop of Lagrangian subspaces in
$\R^{2n}$. Let $\Lambda_n$ denote the Lagrangian Grassmannian of $\R^{2n}$. $H^1(\Lambda_n) \cong \Z$ and there is a distinguished generator
$\mu$ called the Maslov class. The index of the operator is given by applying $\mu$
to the loop of Lagrangian subspaces. Let us denote this index by
$\mu_{\phi}(p,q)$, or by $\mu(p,q)$. Note that this a relative grading in the sense 
that we only know the difference between $p$ and $q$. No absolute grading on the 
critical points is specified.

The space of such $J$-holomorphic discs is at least one-dimensional since we have
a free action by conformal automorphisms of the unit disc $PSL(2,\R)$ on the unit disc. If the
dimension is one, taking the quotient by the free action we expect to get a
zero-dimensional manifold, $\mathcal{M}(p,q)$, as the moduli space of solutions.
A wonderful fact is that there are natural compactifications of spaces of $J$-holomorphic 
maps called Gromov compactifications $\bar{\mathcal{M}}(p,q)$ \cite{MS04}. The codimension one 
component of the boundary of these compactifications, in the best situations, is
\begin{displaymath}
 \del \bar{\mathcal{M}}(p,q) = \coprod_{r \in L\cap L'} \mathcal{M}(p,r) \times
\mathcal{M}(r,q)
\end{displaymath}
This is analogous to finite-dimensional Morse theory where, to compactify the space of gradient trajectories, one adds in trajectories which are broken at an intermediate critical points. The codimension one piece is where we only have one intermediate critical point.
In the case where $\mu(p,q)$ is one, there is no $r$ we can squeeze
in and hence $\mathcal{M}(p,q)$ is compact. Thus, $\mathcal{M}(p,q)$ is a finite
set of points which we can count. Define the chain complex $CF(L,L')$ as the free
$\Z/2\Z$-graded vector space with basis $L \cap L'$ and set 
$$m_1(p) = \sum_{q: \mu(p,q) = 1} n_{pq} q$$
where $n_{pq}$ is the number of points in $\mathcal{M}(p,q)$. 

\begin{prop}
$m_1^2=0$.
\end{prop}

To see this, let $\mu(p,q) =2$ and consider the coefficient of the $q$-term in
$m_1^2(p)$. It is exactly $\sum_r n_{pr} n_{rq}$, which is the number of points in
$\coprod_{r \in L\cap L'} \mathcal{M}(p,r) \times \mathcal{M}(r,q)$.
$\coprod_{r \in L\cap L'} \mathcal{M}(p,r) \times \mathcal{M}(r,q)$ is the
boundary of $\bar{\mathcal{M}}(p,q)$, which is one-dimensional, hence must be an even number. We can
improve the counts to lie in $\Z$, instead of $\Z/2\Z$, when we can coherently orient the moduli
spaces of solutions. This requires a Spin structure \cite{dSi00}.

The cohomology of $CF(L,L')$ is called the (Lagrangian intersection) Floer cohomology 
of $L$ and $L'$ and is denoted by $HF(L,L')$.

What could go wrong with this situation? There are a few problems related to the
Gromov compactification. The compactness relies on bounds on the energy
$\int\phi^*\omega$ of a $J$-holomorphic disc $\phi$. However, in general, we could
have energies of a sequence of $J$-holomorphic discs tending to infinity. To
remedy this, we can include a formal parameter which keeps track of the
symplectic area of these discs. This introduces Novikov rings into the
discussion \cite{HF95}. If we are lucky, as we will be later in this paper, the formal 
series actually converges if we specialize it.

Another, more serious, problem is that if the symplectic form does not vanish on the second homotopy group of $M$, $\pi_2(M)$, 
we would have to include sphere bubbles in the Gromov compactification.
Perturbing the almost complex structure and adding a zeroth order term are no longer
enough to guarantee that the compactified moduli spaces have the proper structure.
This issue was overcome by introducing a new, more general method of perturbation ---
virtual perturbation theory. For the details see \cite{FO99,Ruan98,Sie96,LT98}.

On the other hand, if $\omega$ does not vanish on $\pi_2(M,L)$, we could have disc bubbles which 
generally occupy a (real) codimension one subset of the compactified moduli spaces. Thus, if $\mu(p,q) = 2$, 
the boundary of $\bar{\mathcal{M}}(p,q)$ may now also include other components. Consequently, 
$m_1^2$ may no longer be zero. While this problem cannot always be fixed, there is a nice general 
framework developed in \cite{FOOO} which allows one to address the issue. For a simple and 
illuminating example of what could go wrong in this case, see \cite{Hori02}.

Assume now that our Lagrangian $L$ is sitting inside a Calabi-Yau manifold $X$.
Since $X$ is Calabi-Yau, it possesses a non-vanishing holomorphic volume form $\Omega$. Given a 
Lagrangian submanifold $L$ of $X$, we take a frame 
$\lbrace v_1, \cdots, v_n \rbrace$ for it at a point $x \in X$. This gives a complex number
\begin{displaymath}
 \theta(x) = \frac{\Omega(v_1,\cdots,v_n)^2}{|\Omega(v_1,\cdots,v_n)|^2}
\end{displaymath}
and consequently, a map $\theta: L \ra S^1$ called the phase squared map.
A graded Lagrangian submanifold is a Lagrangian $L$ and a lifting of the
phase map to $\a: L \ra \R$. This induces a grading on the intersection points of graded Lagrangian 
submanifolds. One can check that this grading is compatible with the relative Maslov grading given 
in the Floer complex; see \cite{Sei00}.

A choice of Spin structure, if it exists, and the grading allow one to work, in principle, with $\Z$-coefficients and $\Z$-gradings in Floer cohomology.

There is no reason to stop at $J$-holomorphic discs with two marked points. One can also 
consider $J$-holomorphic discs with $n > 2$ marked points. We outline the construction and refer the inquisitive reader to the references. First we consider 
the universal family of $n$-pointed discs, $\mathcal{A}_n$. By an $n$-pointed disc, we mean a disc 
with $n$ distinguished points on the boundary considered up to equivalence given by M\"obius transformations. 
The marked points $\xi_0,\cdots,\xi_{n-1}$ are ordered cyclically and we label the boundary segment 
between $\xi_i$ and $\xi_{i+1}$ as $C_i$. Given Lagrangians $L_0,\cdots,L_n$ and intersection points 
$p_{ii+1} \in L_i \cap L_{i+1}$ where $i$ is taken modulo $n+1$, we can look for $J$-holomorphic 
maps $\phi$ from an $(n+1)$-pointed disc into our symplectic manifold $M$ so that 
$\phi(C_i) \subset L_i$ and $\phi(\xi_i) = p_{ii+1}$. When $\tilde{p}_{n0} = 2 - n + \sum \tilde{p}_{ii+1}$, 
the space of such maps (when considering all possible $(n+1)$-pointed discs, assuming transversality) is 
zero-dimensional. As before, there are natural compactifications of these moduli spaces which add 
pieces of codimension larger or equal to one. Therefore, in the case where we have chosen the correct intersection 
points, we can count the number of such $J$-holomorphic maps. Denote this count by 
$n_{p_{01}\cdots p_{n-1n} p_{n0}}$ and define 
\begin{displaymath}
 m_i(p_{01},\cdots,p_{n-1n}) = \sum_{q \in L_0 \cap L_n: \tilde{q} = 
2-n + \sum \tilde{p}_{ii+1}} n_{p_{01}\cdots p_{n-1 n} q} q
\end{displaymath}
Recall that $m_1$ satisfied $m_1^2=0$. The total collection of the $m_i$ satisfy the $\Ainf$-relations. 
Let us sketch why this is true. We first consider the case that $M$ is simply a point and look at the 
natural compactifications of these moduli spaces. Any $n$-pointed disc provides a $J$-holomorphic map. 
Thus, our moduli space of such maps is all of $\mathcal{A}_n$. We need to compactify this space. 
To do this we add degenerations of the form

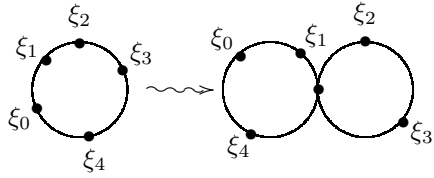
\begin{displaymath}
\xy
 0;/r.30pc/:
 (4,0)*{} = "a";
 (-4,0)*{} = "b";
 (-10,0)*\xycircle(5,5){-};
 {\ar@{~>} "b"; "a"};
 (10,0)*\xycircle(5,5){-};
 (20,0)*\xycircle(5,5){-};
 (15.1,0)*+{\bullet};
 (-5.5,2)*+{\bullet};
 (-14.5,-2)*+{\bullet};
 (-10,4.8)*+{\bullet};
 (-9,-5)*+{\bullet};
 (-13.5,3)*+{\bullet};
 (-16.3,-3)*{\xi_0};
 (-15.5,4.5)*+{\xi_1};
 (-10,7.8)*+{\xi_2};
 (-3.5,3.8)*+{\xi_3};
 (-8.4,-7.4)*+{\xi_4};
 (6.9,3.4)*+{\bullet};
 (8,-4.8)*+{\bullet};
 (13.2,3.7)*+{\bullet};
 (20,5)*+{\bullet};
 (24,-3.5)*+{\bullet};
 (26,-4.5)*+{\xi_3};
 (20,8)*+{\xi_2};
 (5,5.4)*+{\xi_0};
 (14.9,5.7)*+{\xi_1};
 (6.7,-6)*+{\xi_4};
\endxy
\end{displaymath}
as the codimension one boundary strata and work our way down inductively. The codimension one 
boundary strata of $\overline{\mathcal{A}}_n$ are given by 
\begin{displaymath}
 \prod_{n-1>j>2} \mathcal{A}_j \times \mathcal{A}_{n+1-j}
\end{displaymath}
corresponding to the splitting of a single marked disc into two pointed discs, the first with $j$ 
marked points and the second with $n+1-j$ marked points, joined at a marked point. If we view these 
as counting discs, we see that the relations in the codimension one boundary strata are exactly 
the $\Ainf$-relations with vanishing $m_1$. The compactification in the general case works in a very 
similar way. When we generalise from a point, the compactifications include $J$-holomorphic maps with 
two marked points corresponding to the non-vanishing $m_1$. Thus, in the best of circumstances, we can 
form an $\Ainf$-category called the Fukaya category of $M$, $\Fuk(M)$. Its objects are Lagrangian 
submanifolds (with some extra structure) and its morphisms are intersection points between two given 
Lagrangians. The multi-compositions are those outlined above.

We shall see a few variations on the theme of the Fukaya category as we proceed. Most variations 
correspond to restricting one's attention to pieces of a larger Fukaya category. Others correspond to slight 
modifications if we do not actually experience the best of circumstances.

\begin{rmk}
 If one wishes to examine a fuller treatment of the construction of the Fukaya category 
in the case of an exact symplectic manifold, see \cite{SeiDR} and the references therein. For the reader interested in the most general results, see \cite{FOOO}.
\end{rmk}

\section{Homological mirror symmetry for the projective line}
\label{sec:genuszero}

\subsection{B-branes on the projective line}
\label{subsec:BbranesonP1}

We begin by tackling B-branes before A-branes. The name ``branes'' is short for membranes and is a reflection of the subject's physical origin. In string theory, branes are boundary conditions for an open string proprogating through spacetime. The relative simplicity of the algebro-geometric side, or B-side, of mirror symmetry spurs this alphabetic rebellion. After defining the category of B-branes, we seek and find a concrete description of the category in terms of a directed graph, i.e. a quiver.

Let $k$ be a field.

\begin{defn}
The category of B-branes on $\P_k^1$ is the bounded derived category of coherent
sheaves on $\P_k^1$.
\end{defn}

If we are seeking a category that reflects the algebraic geometry of $\P_k^1$,
there is perhaps no better choice. Evidence is given by the following result of Bondal and 
Orlov \cite{BO01}. In passing from $\P_k^1$ to $D^b(\Coh(\P_k^1))$, we lose nothing.

\begin{thm}
Let $X$ and $Y$ be smooth projective varieties and let $\omega_X$ or
$\omega_X^{-1}$ be ample. If $D^b(\Coh(X))$
is equivalent to $D^b(\Coh(Y))$, then $X$ is isomorphic to $Y$.
\end{thm}

In order to get a
better handle on $D^b(\Coh(\P_k^1))$ we seek an alternative description. Consider the structure 
sheaf $\mathcal{O}$ and the twisting sheaf $\mathcal{O}(1)$. Let
$\mathcal{E} = \mathcal{O} \oplus \mathcal{O}(1)$. The coherent sheaf $\mathcal{E}$ ``touches everything'' 
in $D^b(\Coh(\P_k^1))$. More precisely,

\begin{lem}
The smallest triangulated subcategory of $D^b(\Coh(\P_k^1))$ containing $\mathcal{E}$ and closed under taking direct summands is $D^b(\Coh(\P_k^1))$.
\end{lem}

\proof The set of objects $\{\mathcal{O}(n)\}_{n \in \Z}$ generates $D^b(\Coh(\P_k^1))$, since any coherent
sheaf on $\P_k^1$ has a finite resolution by members of this collection (\cite{Hart77}). Consider the following exact sequence 
\begin{displaymath}
 0 \ra \mathcal{O}(n-2) \overset{( y , -x )}{\ra} \mathcal{O}(n-1)^{\oplus 2} \overset{(x,y)^t}{\ra} \mathcal{O}(n) \ra 0
\end{displaymath}
Setting $n=2$ shows that we can resolve $\mathcal{O}(2)$ by $\mathcal{O}(1)$ and
$\mathcal{O}$ and consequently by $\mathcal{E}$. Letting $n=0$ shows that $\mathcal{O}(-1)$ can also 
be resolved (in the other direction) by copies of $\mathcal{E}$. Thus, $\mathcal{O}(2),\mathcal{O}(-1)$ 
lie in the smallest triangulated subcategory of $D^b(\Coh(\P_k^1))$ generated by $\mathcal{E}$ and
closed under taking idempotents. 
Iterating the argument shows inductively that $\mathcal{O}(n)$ also lies in the subcategory 
generated by $\mathcal{E}$, for all $n \in \Z$. Hence, the smallest triangulated
category in $D^b(\Coh(\P_k^1))$ containing $\mathcal{E}$ is $D^b(\Coh(\P_k^1))$ itself. \qed

The next lemma is well known.

\begin{lem}
$\Ext^i(\mathcal{E},\mathcal{E}) = 0$ for $i>0$.
\end{lem}

The only interesting Ext-group is then $\Ext^0(\mathcal{E},\mathcal{E}) = \Hom(\mathcal{E},\mathcal{E})$. 
A simple graphical means of encoding the data of this algebra is the following quiver $Q$.

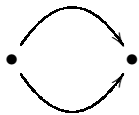
\[
\xy
(-16,0)*+{\bullet}="4";
(0,0)*+{\bullet}="6";
{\ar@/^1.65pc/ "4";"6"};
{\ar@/_1.65pc/ "4";"6"};
\endxy
\]

The path algebra $kQ$ of $Q$ is the algebra generated by oriented paths in $Q$
with multiplication given by concatenation of paths. Note that the trivial paths at each of the 
vertices are idempotents.

\begin{lem}
$\Hom(\mathcal{E},\mathcal{E}) \cong kQ$.
\end{lem}

\proof $\mathcal{O}$ corresponds to the first node and $\mathcal{O}(1)$ to the
second. The two arrows are the morphisms $x,y: \mathcal{O} \ra \mathcal{O}(1)$. \qed

From what we have just seen, the following should not be too surprising.

\begin{prop}
 $D^b(\Coh(\P_k^1))$ is triangle equivalent to $D^b(\mbox{$kQ$-$\mod$})$.
\end{prop}

The proof of this statement requires a little sophistication. Recall that, to get $D^b(\Coh(\P_k^1))$, 
we replaced all our bounded complexes of coherent sheaves with choices of injective resolutions and then took the homotopy category. Let 
us take an injective resolution of $\mathcal{E}$ and denote it by $I_{\mathcal{E}}$ and consider 
the category $I(\P_k^1)$ formed by injective complexes with bounded coherent cohomology. We then consider the functor 
$\Hom_{I(\P_k^1)}(I_{\mathcal{E}},-)$ from $I(\P_k^1)$ to $Ch(k)$. In fact, the image of an object 
$J$ in $Ch(k)$ has the structure of a left dg-module (and hence $\Ainf$-module) over the dg-algebra 
$A = \Hom_{I(\P_k^1)}(I_{\mathcal{E}},I_{\mathcal{E}})$ given by precomposition. Thus, we have a 
functor from $I(\P_k^1)$ to $A$-\mbox{Mod}. This induces a functor from $I(\P_k^1)$ to $D(A)$. The question 
then becomes: does it descend to a functor from $D^b(\Coh(P_k^1))$ to $D(A^{op})$? In other words, are 
all quasi-isomorphisms in $I(\P_k^1)$ inverted by the functor to $D(A^{op})$? The answer is yes. 
Since any quasi-isomorphism has an inverse up to homotopy in $I(\P_k^1)$, we just need to check that 
our functor kills all null-homotopic maps. This is equivalent to checking that the functor kills all 
null-homotopic complexes, i.e. those complexes for which the identity map is null-homotopic. But the 
contracting homotopy for the identity map of the complex, $N$, provides the contracting homotopy for 
any chain map in $\Hom_{I(\P_k^1)}(I_{\mathcal{E}},N)$. Thus, the functor descends. Let us denote 
the resulting functor by $\RHom(\mathcal{E},-)$. The dg-algebra $A$ is quasi-isomorphic to $kQ$. 
Applying Proposition \ref{prop:quasiequivcat}, we can assume that $\RHom(\mathcal{E},-)$ maps from $D^b(\Coh(\P^1_k))$ to $D(kQ^{op})$. $\RHom(\mathcal{E},-)$ is full and faithful on 
the smallest triangulated category containing $\mathcal{E}$ and closed under direct summands and hence is fully faithful. It is 
also essentially surjective onto $D^{\pi}(kQ^{op})$ since both have to be the smallest 
triangulated categories containing $kQ$ and closed under direct summands. Now, as is well known, $kQ$ has finite 
global dimension thus, $D^{\pi}(kQ^{op}) \cong D^b(\mbox{$kQ$-{mod}})$; see \cite{Rou03}.

Thus, the structure of the seemingly complicated $D^b(\Coh(\P_k^1))$ is controlled
by the relatively simple algebra $kQ$. We shall see the quiver $Q$ again shortly.

\subsection{A-branes on the mirror to the projective line}
\label{sec:AbranesMirrorP1}

The mirror partner for $\P^1_k$ is something a little more exotic than a variety. 
We consider the function $W: \C^{\times} \ra \C$ given by $W(z) = z + q/z$
with $q \in \C^{\times}$. Such a pair $(\C^{\times},W)$ is called a Landau-Ginzburg model, or LG-model. 
Since the value of $q$ does not affect the category of A-branes, we set $q=1$. When dealing with LG-models, a useful moral to keep in mind is 
that we should study the geometry of the critical locus and how it relates to the ambient 
space. Now, on the A-side, we are interested in symplectic geometry, in particular Lagrangians. First, 
we need a symplectic form. We take it to be 
$\displaystyle{\imath \frac{dz \wedge d\bar{z}}{z\bar{z}}}$.

We could just take Lagrangian submanifolds of $\C^{\times}$ as our A-branes but that would not reflect 
the presence of $W$. Instead, we will look at non-closed Lagrangian submanifolds of $\C^{\times}$ whose 
boundaries lie on $W^{-1}(z)$ for some generic choice $z$. We could proceed as in \cite{Abou07} and 
investigate all such Lagrangians with some additional requirements. This approach seems much more 
natural, but for simplicity's sake we shall follow the ideas of Seidel in \cite{Sei001} or \cite{AKO04} 
and only consider a few special Lagrangians. $W$ has non-degenerate singularities, i.e., for each $p$ such that $W'(p)=0$ we have $W''(p) \not = 0$. Each such $p$ is called a critical point; $W(p)$ is 
called the critical value. There is a natural connection on the complement of the singular fibers 
coming in this situation. Given a vector field $X$ on $\C$ we can find a horizontal lift $\tilde{X}$ 
by using the symplectic form to split the tangent space of the domain of $W$ into the tangent space of the fiber and its 
symplectic orthogonal. We take $\tilde{X}$ to lie in this orthogonal and project down to $X$. With 
a connection, we can parallel transport vectors. Thanks to our choice of connection, parallel 
transport preserves the symplectic form. This gives our Lagrangians. Fix a non-critical 
point $q$. If we take a path $\gamma: [0,1] \ra \C$ so that $\gamma(0) = p$ and $\gamma(1)=q$, we can 
look at the set
\begin{displaymath}
 \lbrace x \in W^{-1}(q) : \lim_{t \ra 0} T_{\gamma|_{[t,1]}}(x) = q \rbrace \cup \lbrace q \rbrace
\end{displaymath}
where $T_{\gamma|_{[t,1]}}(x)$ is the symplectomorphism coming from parallel transport along 
$\gamma|_{[t,1]}$. The resulting set is just an interval (or a one-dimensional disc). It is a Lagrangian submanifold with boundary on $W^{-1}(q)$. 
Stated in this language, the fact that the stable manifold is an $n$-dimensional Lagrangian disc
holds in greater generality; see \cite{Sei001}. This Lagrangian is called the vanishing thimble associated 
with $p$. The boundary, which is an $(n-1)$-dimensional sphere, is called the vanishing cycle associated 
with $p$. Since our example is low-dimensional, we can actually draw the vanishing thimbles associated 
to the two critical points for $z + 1/z$.

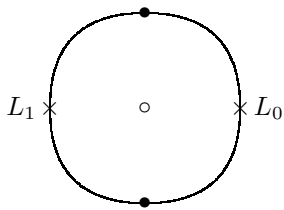
\begin{displaymath}
 \xy
  0;/r.30pc/:
  (0,0)*+{\circ};
  (-10,0)*+{\times};
  (10,0)*+{\times};
  (0,10)*+{\bullet};
  (0,-10)*+{\bullet};
  (0,10)*{}="A";
  (10,0)*{}="B";
  "A"; "B" **\crv{(10,10)};
  (-10,0)*{}="C";
  "A"; "C" **\crv{(-10,10)};
  (0,-10)*{}="D";
  (-10,0)*{}="E";
  "D"; "E" **\crv{(-10,-10)};
  (10,0)*{}="F";
  "D"; "F" **\crv{(10,-10)};
  (13,0)*+{L_0};
  (-13,0)*+{L_1};
 \endxy
\end{displaymath}

We can now define our category of A-branes associated with the LG-model.

Our two critical points give us two Lagrangian thimbles which intersect in two points. 
Order the Lagrangian thimbles $L_0$ and $L_1$. 

\begin{defn}
The category of vanishing cycles, $FS(W)$ for the Landau-Ginzburg model $W$ is an $\Ainf$-category 
whose objects are collection of vanishing thimbles $L_i$,  and whose morphisms are
\begin{displaymath}
 \Hom(L_i,L_j) = \begin{cases} CF(\del L_i,\del L_j) & \text{ if $i < j$} \\ k
\id_{L_i} & \text{ if $i=j$} \\ 0 & \text{ otherwise} \end{cases}
\end{displaymath}
We use the usual $\Ainf$-structure on the vanishing cycles of the regular fiber $W=0$.
The category of A-branes is the bounded derived idempotent completion of the category of vanishing 
cycles, $D^{\pi}FS(W)$.
\end{defn}

One should first notice that the $\Ainf$-relations respect the (strict) ordering, in the sense that 
any degeneration of an ordered disc splits into two ordered discs. Thus, the $\Ainf$-structure on the 
Fukaya category gives $FS(W)$ its $\Ainf$ structure.

One should also note that, a priori, $D^{\pi}\FS(W)$ depends on the choice of paths. Different choices 
yield different Lagrangians. But, in fact, $D^{\pi}FS(W)$ is independent up to equivalence 
of the choice of paths \cite{Sei002}.

\begin{defn}
 The category of A-branes for the LG-model $W$ is the category $D^{\pi}FS(W)$.
\end{defn}

Now, let us take a closer look at the category of A-branes for $W(z) = z+1/z$. We have two objects, 
two morphisms from the first object to the second, and the identity morphisms. Thus, it almost a 
tautology that the morphism algebra of $FS(W)$ is the path algebra of our quiver $Q$.  
$D^{\pi}(FS(W))$ is simply the bounded derived category of finite-dimensional representations of $Q$. 
Indeed, the category of B-branes on $\P^1_k$ is equivalent to the category of A-branes on $W$. 

\begin{rmk}
 The simplicity of this example hides several important details. One of these details is that 
$\mathcal{O},\cdots,\mathcal{O}(n)$ is the wrong exceptional collection to look at. One should 
instead consider $\Omega^n(n),\cdots,\Omega(1),\mathcal{O}$. The case of the mirror to $\P^2_k$ 
provides another illuminating example of homological mirror symmetry; see \cite{AKO04}.
\end{rmk}

\subsection{B-branes on the mirror to the projective line}
\label{sec:BbranesonMirrorP1}

With half of the homological mirror symmetry correspondence done, we now turn
our attention to the other half. The naive guess for B-branes on the Landau-Ginzburg theory 
$W: \C^{\times} \ra \C$ would simply be coherent sheaves on
$\C^{\times}$ or equivalently modules over $\C[z,z^{-1}]$. What we should remember here is that, 
as in the case of A-branes on the Landau-Ginzburg model, we want to take the potential $W$ into account. We must measure the singularities of $W$ using complex geometry. This reflects 
the similar considerations on the A-side where the singular fibers of $W$ gave rise to the vanishing 
thimbles. How do we measure singularities? By a classical result of Serre, a Noetherian commutative algebra $A$ is 
regular if and only if every module over $A$ has a projective resolution of uniformly bounded length. 
For a general $A$, modules with such bounded resolutions form a triangulated
subcategory of $D^b(\mbox{{mod}-$A$})$. Let us denote this subcategory by $\Perf(A)$. Then
a measure of the singularity of $A$ is the quotient $D^b(\mbox{{mod}-$A$})/\Perf(A)$. This is often called the 
stable category of $A$ or the category of singularities of $A$ and it is indeed the desired definition 
of B-branes on $W$.

\begin{defn}
The category of B-branes of $W: Y \ra \C$ is
\begin{displaymath}
 D_{\rm{sing}}(W) = \prod_{\lambda \in \C} D^b(\Coh(Y_{\lambda}))/\Perf(Y_{\lambda})
\end{displaymath}
where $\Perf(Y_{\lambda})$ is the full triangulated subcategory of complexes of coherent sheaves
admitting a bounded locally free resolution.
\end{defn}

Note that, since we are working over $\C$, there are only a finite number of singular values 
and hence the product is finite. Of course, for this definition to be useful, one needs to know 
how to take the quotient of $D^b(\Coh(Y_{\lambda}))$ by $\Perf(Y_{\lambda})$. We will not cover this, because there is a more computationally accessible version of the category of B-branes on $W$ --- the category of matrix factorizations of $W$, $MF(W)$.

Let $A$ be an algebra and $f \in A$. A matrix factorization of $f$ is a periodic
complex 
\begin{displaymath}
 \cdots \ra P_0 \overset{d_P^0}{\ra} P_1 \overset{d_P^1}{\ra} P_0
\overset{d_P^0}{\ra} P_1 \ra \cdots
\end{displaymath}
with $P_0,P_1$ projective modules over $A$ and $d_P^0d_P^1 = d^1_Pd_P^0 = f\id$. Notice
that we could repackage the data by considering $P= P_0 \oplus P_1$ with $d_P = 
\begin{pmatrix}  0 & d_P^0 \\ d_P^1 & 0 \end{pmatrix}$ so that  $d^2 = f\id$. Complexes of this form will be the objects in the category of matrix factorizations of $f$.
Given two matrix factorizations $P^{\bullet}$ and $Q^{\bullet}$, we define the set of
morphisms  $\Hom_{CF(f)}(P^{\bullet},Q^{\bullet})$ to be $\Hom_A(P,Q)$.  It is 
$\Z/2\Z$-graded, with the degree zero piece being given by $\Hom_A(P_0,Q_0) \oplus
\Hom_A(P_1,Q_1)$ and the degree one piece being given by $\Hom_A(P_0,Q_1) \oplus \Hom_A(P_1,Q_0)$. One can view this as a periodic version of the ordinary Hom-complex. One nice 
feature is that while $P^{\bullet}$ and $Q^{\bullet}$ are not chain complexes, 
$\Hom(P^{\bullet},Q^{\bullet})$ is. The action of the differential $d_{P,Q}$ is simply super-commutation, 
namely $d_{P,Q} \phi = d_Q \phi - (-1)^{\tilde{\phi}}\phi d_P.$ The proof of the following result is left as a simple exercise to the reader.

\begin{lem}
The composition $d^2_{P,Q}$ is $0$. Thus it makes sense to speak of the cohomology of $\Hom(P^{\bullet},Q^{\bullet})$. 
\end{lem}

We also have a shift functor $[1]$ which sends $P_i$ to $P_{i+1 \mod 2}$ and
$d_i$ to $-d_{i+1 \mod 2}$ and which acts on morphisms in the standard way. Notice that $[2] \cong \Id$. 
We can form cones over morphisms $\phi: P^{\bullet} \ra Q^{\bullet}$ by 
taking $C(\phi) =  P[1] \oplus Q$ 
with factorizations $d_{C(\phi)}=\begin{pmatrix} d_{P[1]} & \phi \\ 0 & d_Q \end{pmatrix}$.

\begin{defn}
The category of matrix factorizations of $f$, $MF(f)$, is the category whose
objects are matrix factorizations and whose morphisms are
\begin{displaymath}
 \Hom_{MF(f)}(P^{\bullet}, Q^{\bullet}) = H^0(\Hom_{CF(f)}(P^{\bullet},Q^{\bullet}))
\end{displaymath}
\end{defn}

\begin{prop}
$MF(f)$ is a triangulated category where the shift is as above and triangles are
of the form
\begin{center}
\leavevmode
\begin{xy}
 (-10,10)*+{A}="a"; (10,10)*+{B}="b"; (0,-5)*+{C(\phi)}="c"; {\ar@{->}^{\phi} "a";"b"}; {\ar@{->} "b";"c"}; {\ar@{->}^{[1]} "c";"a"}
\end{xy}
\end{center}
\end{prop}

This proposition is proved in \cite{Orl04}, where we find the following very
useful theorem.

\begin{thm}
Consider a Landau-Ginzburg model $W: Y \ra \C$ with $Y$ affine. Then $D_{\rm{sing}}(W)
= \prod_{\lambda \in \C} MF(W - \lambda) =: MF(W)$ as triangulated
categories.
\label{thm:Orlov}
\end{thm}

\begin{rmk}
This is an extension of the proof by Eisenbud that the category of matrix
factorizations is equivalent to the category of Cohen-Macaulay modules
\cite{Eis80}.
\end{rmk}

Here is an elucidating example.

\begin{eg}
Let $A = k[x]/(x^2)$. Then the minimal free resolution of $A/(x)$ is 
$$ \cdots A \overset{x}{\ra} A \overset{x}{\ra} A \ra A/(x) \ra 0$$
Lifting this resolution to $k[x]$ we get a matrix factorization of $x^2$.
\end{eg}

In fact, our situation is quite similar to the previous example. In the case of
the mirror to $\P_k^1$, $W = z + q/z$. The critical points of $W$ are $\pm \sqrt{q}$. The critical values are $z = \pm 2\sqrt{q}$ and $W \pm 2\sqrt{q} = (1/z)(z\pm\sqrt{q})^2$. Let us take
the matrix factorization, $F$, for $W-2\sqrt{q}$, given by $P_0 = P_1 =
\C[z,z^{-1}]$ and $d_P^0 =(1/z)(z-\sqrt{q}), d_P^1 = (z - \sqrt{q})$. First,
let us compute the cohomology algebra $\Hom_{MF(W-2\sqrt{q})}(F,F)$.

\begin{prop}
$\Hom_{MF(W-2\sqrt{q})}(F,F)$ is isomorphic to the Clifford algebra on a single
generator $\phi$ with $\phi \cdot \phi = -(\del^2_z W)(\sqrt{q})$.
\end{prop}

\proof A general morphism $\psi \in \Hom_{CF(W-2\sqrt{q})}(F,F)$ is of the form
$\psi = \begin{pmatrix} A & B \\ C & D \end{pmatrix}$ and 
$$d_{F,F} \psi=$$
$$\begin{pmatrix}0 & (1/z)(z-\sqrt{q}) \\ (z-\sqrt{q}) &
0\end{pmatrix} \begin{pmatrix} A & B \\ C & D\end{pmatrix} - \begin{pmatrix} A &
-B \\ -C & D \end{pmatrix} \begin{pmatrix} 0 & (1/z)(z-\sqrt{q}) \\ (z-\sqrt{q}) &
0 \end{pmatrix}$$
$$=\begin{pmatrix} B/z+C & D-A \\ B/z+C & A-D \end{pmatrix} \begin{pmatrix}
z-\sqrt{q} & 0 \\ 0 & z-\sqrt{q} \end{pmatrix}$$
Thus for $\psi$ to lie in the kernel we must have $A=D$ and $B=-Cz$. Over
$\C[z,z^{-1}]$, the cohomology is generated by the identity
and by $\phi = \begin{pmatrix} 0 & 1/z \\ -1 & 0 \end{pmatrix}$. Note that
$\phi^2 = \begin{pmatrix} -1/z & 0 \\ 0 & -1/z \end{pmatrix}$. Notice that
$d_F^2 = W-\lambda$ implies that
$d_F (\del_z d_F) + (\del_z d_F) d_F = \del_z W$. Hence $\del_z W \id$ is
trivial in cohomology. Thus, the cohomology algebra
becomes a module over the Jacobian ring $\C[z,z^{-1}]/(\del_z W)$, which in
particular is finite dimensional over $\C$. Now $\phi^2 = (-1/z) \id$.
The image of $1/z$ in the Jacobian ring is $1/\sqrt{q}$ and hence $\phi^2 =
-(\del_z^2 W)(\sqrt{q})$. \qed

A similar result holds at $z = - \sqrt{q}$. 

Usually, the cohomology algebra is not enough data. We need to know something about 
the underlying $\Ainf$-structures. However, in this case, we actually have all the data we need.
Hochschild cohomology of a graded associative algebra $A$ can be computed as $\Ext^*_{A^e}(A,A)$.

\begin{prop}
 A Clifford algebra $C$, considered as a $\Z/2\Z$-graded algebra, is projective over $C^e$.
\end{prop}

We shall not recall the proof here, but only refer the reader to \cite{Bass67} for satisfaction.

Appealing to the previous result and lemma \ref{lem:HochtrivAinf}, we can conclude that the dg-algebra $\Hom_{CF(W-2\sqrt{q})}(F,F)$ is quasi-isomorphic to its cohomology.
Again a similar result holds at $z = - \sqrt{q}$.

Finally, we can appeal to Orlov's result, theorem \ref{thm:Orlov}, to conclude that the zero object and $F$ generate 
$MF(W-2\sqrt{q})$, since the fiber over the singular points is isomorphic to $k[x]/(x^2)$. We can 
conclude that $MF(W-2\sqrt{q})$ is equivalent to the bounded derived category of modules over 
the Clifford algebra $\Hom_{MF(W-2\sqrt{q}}(F,F)$ via an argument analogous to that in subsection 
\ref{subsec:BbranesonP1}. Thus, the category of B-branes for $z+q/z$ is equivalent to two copies 
of the bounded derived category of modules over the Clifford algebras 
$k\left<x\right>/\left<x^2=\pm1/\sqrt{q}\right>$.

\subsection{A-branes on the projective line, a.k.a. the sphere}
\label{subsec:A-branesonP1}

Equip $\P_{\C}^1$ with the Fubini-Study metric and its associated K\"ahler form
$\omega$. In a Riemann surface, any curve is a Lagrangian. We shall
choose to deal only with embedded and closed Lagrangians, i.e. simple closed
curves in $\P_{\C}^1$. Any such curve divides the sphere into two, perhaps unequal (with respect 
to the Fubini-Study area form), halves. We first reproduce an illustrative calculation from 
\cite{Hori02}. Assume we have two Lagrangians $\gamma$ and $\gamma'$ that intersect in exactly two 
points. Let us denote the regions as in the image below. 

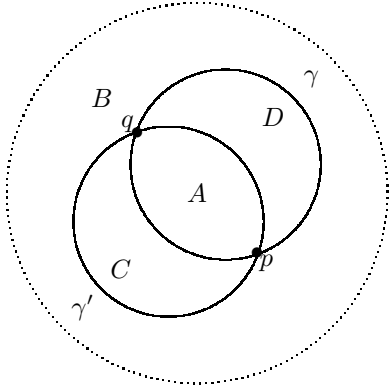
\begin{displaymath}
 \xy
 0;/r.30pc/:
 (-3,-3)*\xycircle(10,10){-};
 (3,3)*\xycircle(10,10){-};
 (0,0)*+{A};
 (-10,10)*+{B};
 (-8,-8)*+{C};
 (8,8)*+{D};
 (-6.3,6.3)*+{\bullet};
 (-7.3,7.3)*+{q};
 (6.3,-6.3)*+{\bullet};
 (7.3,-7.3)*+{p};
 (12,12)*+{\gamma};
 (-12,-12)*+{\gamma'};
 (0,0)*\xycircle(20,20){.};
 \endxy
\end{displaymath}

Let $u$ be a formal parameter that we will use to keep track of the symplectic area
of the discs of interest. We shall use the standard complex structure on the sphere. 
Before embarking on the computation, we must make sure that we have the signs associated to our 
Floer differential correct and we must make sure that we have the proper gradings on the intersection 
points. A thorough discussion of a coherent sign convention would take us off track, so we shall carry 
out our computations with $\Z/2\Z$ coefficients. This will initially eliminate considerations of the 
line bundles supported along the Lagrangians. The issue of grading, however, cannot be avoided. To simplify 
the situation, here we shall assume that one of our Lagrangians is a small Hamiltonian deformation 
associated to a Morse function. In this case, the intersection points are identified with the critical 
points of the Morse function and we can use the associated Morse index to grade them.

Now we compute the Floer differential.
\begin{displaymath}
 m_1(q) = (u^{\omega(C)}+u^{\omega(D)})p
\end{displaymath}
\begin{displaymath}
 m_1(p) = (u^{\omega(A)}+u^{\omega(B)})q
\end{displaymath}
Then,
\begin{displaymath}
 m_1^2(p) = (u^{\omega(C)}+u^{\omega(D)})(u^{\omega(A)}+u^{\omega(B)}) p =
\end{displaymath}
\begin{displaymath}
 u^{\omega(C)+\omega(B)}(u^{\omega(\Int(\gamma))}+u^{\omega(\Int(\gamma'))})
(u^{\omega(\Int(\gamma'))}+u^{\omega(\Out(\gamma)})p
\end{displaymath}

Since $\gamma'$ is a Hamiltonian deformation of $\gamma$ the areas of the interiors
of $\gamma$ and $\gamma'$ are equal. Hence, $m_1^2(p) = 0$ and similarly
$m_1^2(q)=0$. Thus, the Floer cohomology between $\gamma$ and $\gamma'$ is well-defined. Note that 
the equality of the interior areas implies that $\omega(A)=\omega(B)$.
Thus, $m_1(p) = 0$. Now consider $m_1(q)$.
\begin{displaymath}
 m_1(q) = (u^{\omega(C)}+u^{\omega(D)})p = u^{\omega(C)}(1+u^{\omega(D)-\omega(C)})p =
\end{displaymath}
\begin{displaymath}
 u^{\omega(C)}(1+u^{\omega(\Out(\gamma))-\omega(\Int(\gamma'))})p =
\end{displaymath}
\begin{displaymath}
 u^{\omega(C)}(1+u^{\omega(\Out(\gamma))-\omega(\Int(\gamma))})p
\end{displaymath}

Thus, unless $\gamma$ cuts the sphere in half, $m_1(q) \not =0$ and the Floer cohomology of 
$\gamma$ with itself vanishes. In the cohomological category of the Fukaya category 
of $\P_{\C}^1$, $\gamma$ is isomorphic to the zero object if and only if $\gamma$ separates 
the sphere into two unequal halves.

The most obvious choice of a $\gamma$ which is cohomologically non-zero is a great circle. It is 
a classical result of Poincar\'e that any length minimizing area-bisecting curve in the sphere must 
be a great circle. The good news is that any
area-bisecting curve is Hamiltonian isotopic to a great circle \cite{Oh90}.
Thus, up to isomorphism the only nontrivial objects in the Fukaya category for
$\P_{\C}^1$ are the great circles. Given two great circles, we can take the poles
of the sphere to be their intersection points and take the height function $h$ on
the sphere. The associated Hamiltonian vector field is simply $J\nabla h$ where
$J$ is the complex structure on $\P_{\C}^1$. Thus, the associated Hamiltonian flow
will simply rotate the sphere along the given axis and eventually take one great
circle into another.

From our brusque analysis, we have concluded that with $\Z/2\Z$-coefficients there
is only one (cohomologically) non-zero object in the Fukaya category --- a choice
of a great circle. If we keep trivial flat line bundles, this continues to hold
with $\Z$ coefficients. We shall switch to $\C$ coefficients and outline the
computation showing that the holonomy of our flat line bundle must be $\pm 1$ in order 
for the A-brane to be non-zero. The computation showing that $m_1^2=0$ remains essentially 
the same. In the presence of a non-trivial flat line bundle, we get the following equation
\begin{displaymath}
 m_1(q,v) = u^{\omega(C)}(v-e^{\int_{\gamma'}A'}ve^{\int_{\gamma}A} u^{\omega(\Out(\gamma))
-\omega(\Int(\gamma))})p =
\end{displaymath}
\begin{displaymath}
 u^{\omega(C)}(v-(\Hol_{\gamma}(A))^2 v u^{\omega(\Out(\gamma))
 -\omega(\Int(\gamma))}p
\end{displaymath}
Thus, the holonomy of the flat line bundle must be $\pm 1$ and as before $\gamma$
must bisect $\P_{\C}^1$. Consequently, we have two non-trivial objects in the Fukaya category --- 
a great circle with either trivial line bundle or a flat line bundle with holonomy $-1$. To determine 
the full structure of the category we must compute the morphism spaces. A computation similar to 
the above shows that $\Hom((\gamma,L_1),(\gamma,L_{-1}))$ has zero cohomology. Thus, there are two 
non-zero objects which are orthogonal to each other, at least cohomologically. The structure of 
the category of A-branes on $\P_{\C}^1$ now begins to resemble the structure of the category of 
B-branes on the mirror to $\P_{\C}^1$.

It is important to keep in mind that the above discussion sidestepped many of the
crucial points in the construction of the Fukaya category, such as whether it
honestly factors through Hamiltonian isotopy, orientations of the relevant
moduli spaces, transversality, checking the $\Ainf$-structure, etc. Many of
these issues are resolved.
For instance, the issues of Hamiltonian isotopy and transversality are taken
care of in this case; see \cite{CO06}.

Now that we have seen that there are really only two non-zero objects of the
Fukaya category of $\P_{\C}^1$ and that the only non-zero morphism sets are their
endomorphisms, we will present another way of computing $HF(L,L)$ which was
originally conceived in \cite{FOOO}. It is a natural adaptation of Bott's
generalisation of Morse theory to allow non-isolated critical points. Naturally,
it is called Bott-Morse-Floer cohomology. As mentioned above, bubbling 
holomorphic discs can occur in codimension one strata of the moduli space of
holomorphic discs. The resolution of this problem was the central motivation of
\cite{FOOO}. The presence of nontrivial discs with boundary on a given
Lagrangian should deform the singular cohomology algebra structure in much the
same way that the presence of nontrivial holomorphic spheres deforms the singular
cohomology into the (small) quantum cohomology. Consider the singular chain
complex of our great circle with its usual intersection product. The
intersection product is not defined on all singular chains so we must restrict
the allowed class of composable objects to those that are transversal. As the
reader might appreciate now, the issue of transversality, or more precisely the
lack thereof, crops up often in these cases. One means to sidestep the issue is
the definition of a pre-category \cite{KS01,Abou07}. While we shall not be so
formal, it is important to keep this point in mind. We now outline the
computation of Bott-Morse-Floer cohomology obtained in \cite{Cho03,Cho05}.

Let $\mathcal{M}_k(\beta)$ be the moduli space of stable maps $f$ from
$(k+1)$-marked discs to $\P_{\C}^1$ with $f_*[\del D] = \beta \in \pi_2(\P_{\C}^1,L)$. 
Here we have cyclically ordered the marked points. The markings provide us with
evaluation morphisms $ev_i: \mathcal{M}(\beta) \ra L$. Now, we deform the
intersection product as follows: $$m_2(S,T) = S \cap T + \sum_{\beta \in
\pi_2(\P,L)} ev_0(ev_1^*S \cap ev_2^*T) u^{\omega(\beta)}.$$ As the reader might
have guessed by now, the deformed composition is no longer associative. We must
add in compositions of all orders to obtain an $\Ainf$-algebra. In fact, in this
case, we must add in an $m_0$ also. We set $$m_0 = \sum_{\beta \in \pi_2(\P,L)}
ev_0(\mathcal{M}_0(\beta))u^{\omega(\beta)}$$ and
\begin{displaymath}
 m_k(S_1,\cdots,S_k) = \sum_{\beta \in \pi_2(\P^1_{\C},L)} ev_0(\cap ev_i^*S_i)
u^{\omega(\beta)}
\end{displaymath}
In fact we only get a $\Z/2\Z$-grading on the resulting $\Ainf$-algebra. If we
added in another formal parameter counting the Maslov index of $\beta$, we could
rectify this. 

The proper way to view expressions such as $ev_0(\cap ev_i^*S_i)$ is as currents,
i.e. distribution-valued differential forms. Thus, the image, $ev_0(\cap
ev_i^*S_i)$, will represent the same current as its closure. 

\begin{prop}
$m_k(x_1,\cdots,[L],\cdots,x_{k-1}) = 0$ for $k \not = 2$ and $m_2([L],x) =
(-1)^{1-\tilde{x}} m_2(x,[L]) = x$.
\end{prop}

\proof It is clear that the fundamental class is closed since any quantum
correction gives a chain with too large a degree. Since any class meets
$L$, the class $$ev_0(ev_1^* L \cap ev_i^*S_i)$$ lies inside $ev_0(\cap ev_i^*S_i)$. If $k>2$,
the dimension of $\mathcal{M}_{k-1}(\beta)$ is one less than the dimension of
$\mathcal{M}_k(\beta)$ since we have lost a marked point and consequently a
degree of freedom. Thus, the actual dimension of the current is smaller than the
expected dimension and the current vanishes. The same reasoning tells us that in the
case $k=2$ any quantum corrections vanish leaving only the standard cap product.
The sign arises from the choice of convention in \cite{Cho05}. \qed

Notice that $m_0$ is a multiple of the fundamental class. As we remarked previously, 
we still have an honest $\Ainf$-algebra structure on cohomology. Since $L$ is a strict unit, 
the only compositions
that matter are those involving chains in the class of points. Now we check the
previous calculation concerning the (Floer-cohomological) non-triviality of a
given Lagrangian.

\begin{prop}
Unless the Lagrangian divides $\P_{\C}^1$ in half, the Floer cohomology vanishes.
\end{prop}

\proof $m_1(p) = [L](u^{\omega(\beta_1)}-u^{\omega(\beta_2)})$. The reversal in
the sign comes from the different orientations induced on the boundary from each
attaching disc. With a single marked point we trace out the whole of $L$. Now, if
the symplectic area associated to $\beta_1$ is not equal to that $\beta_2$, we
see that $m_1(p)$ is a nonzero multiple of the $[L]$. Note that, for any
one-chain $S$, $m_1(S)$ contains a piece corresponding to the boundary of $S$ as
a chain. The rest is a singular chain of dimension one. Hence, for a singular
chain of dimension one to be in the kernel of $m_1$, the chain needs to be the
fundamental class or some multiple of it. Thus, if the Lagrangian
does not cut $\P_{\C}^1$ in half, the Floer cohomology is zero. \qed

The argument in the previous proposition also shows that, if we take the equator as Lagrangian,
the Floer cohomology is isomorphic to the standard cohomology as a module over
our ring $\C[u]$.  If we equip the Lagrangian with a flat
line bundle $\mathcal{L}$, we incorporate the holonomy into the
$\Ainf$-structure as follows
\begin{displaymath}
 m_k(S_1,\cdots,S_k) = \sum_{\beta \in \pi_2(\P^1_{\C},L)} ev_0(\cap ev_i^*S_i)
u^{\omega(\beta)}\hol_{\del \beta}(\mathcal{L})
\end{displaymath}
As before, we have reduced to the two A-branes of interest: the equator equipped
with the trivial line bundle and the equator equipped with the non-trivial flat
line bundle. In these two cases, the Floer cohomology is non-vanishing and
isomorphic as a module to the usual cohomology. The richness lies in the deformation of the
standard algebra structure. In our case, the usual exterior
algebra structure of the cohomology of a circle is deformed to a Clifford
algebra $k\left<x\right>/\left<x^2 = \pm u^{\omega(\P^1_{\C})/2}\right>$. 

\begin{prop} 
Let $p,q$ be distinct points thought of as transversal chains. Then 
$$m_2(p,q) = [L]e^{\omega(\beta)}\hol_{\del \beta}(\mathcal{L})$$ 
\end{prop} 
 
\proof When we consider a disc with three marked points two of which must lie on 
$p$ and $q$, we can either orient the disc so that the zeroth marked point is 
before or after $p$. The two orientations are realized by the two discs 
attaching to $L$. Thus, we trace out all of $L$. The subtlety is in the signs. The 
moduli spaces of the discs all have natural orientations from auxiliary choices. 
In this case the induced orientations match up and we get all of the fundamental 
class. \qed 

Stopping and taking stock of the situation, we see that on the A-side and B-side 
we have two special objects whose $\Z/2\Z$-graded endomorphism algebras are 
Clifford algebras.

If we specialize our formal parameter $u$ to $e$, as is physically motivated, we end up 
in the same situation as before. Any $\Ainf$-structure on $HF(L,L)$ is necessarily 
isomorphic to the trivial one. When we pass to $D^{\pi}(\Fuk(\P^1_{\C}))$ we will get 
a sum of the bounded derived category of modules over the two Clifford algebras 
$k\left<x\right>/\left<x^2 = \pm e^{\omega(\P^1_{\C})/2}\right>$.

\begin{quoteprop}
The category of A-branes on $\P_C^1$ is equivalent to the category of B-branes on
$W:\C^{\times} \ra \C, W(z) = z + q/z$ where $q = e^{-\omega(\P^1_{\C})}$. Mathematically, the 
derived idempotent-completed Fukaya category of $\P_k^1$ is equivalent to the triangulated 
category of matrix factorizations of its mirror.
\label{prop:AtoB}
\end{quoteprop}

\begin{rmk} We can interpret this correspondence of A-branes on $\P^1$ and 
B-branes on $W$ as an instance of T-duality. T here is for the torus. $\P^1$ admits a 
torus fibration over the interval given by the momentum map of the 
$U(1)$-rotation along an axis. In T-duality, we replace these tori with their 
dual tori, i.e. $\Hom(T,U(1))$. Often, we think about the torus as having some metric 
structure which is a choice of length $R$. The T-dual torus will then have 
length $1/R$. Thus, if we T-dualize the momentum map torus fibration for $\P^1$, 
the two ends open up and we get the algebraic torus $\C^{\times}$. The two 
non-trivial A-branes map to $1$ and $-1$ as determined by their holonomies. These are exactly the critical points of the superpotential on $\C^{\times}$.  
\end{rmk}

\begin{rmk}
 As mentioned before, there are a few details missing from a complete proof of proposition 
\ref{prop:AtoB}, such as the complete definition of $\Fuk(\P^1_{\C})$! A rigorous reformulation 
of the previous proposition would just require an assumption that a $\Fuk(\P^1_{\C})$ exists 
such that the computations carried out above are valid.
\end{rmk}

\section{Homological mirror symmetry for elliptic curves}
\label{sec:genusone}

The name ``mirror symmetry'' implies that two objects related by the correspondence should be 
similar enough for one to be a reflection of the other. In the previous section we saw that 
$\P_k^1$ had as mirror not another algebraic variety, not even another topological space, but 
a function. To the untrained eye, these two objects appear quite dissimilar. For something a little 
more symmetric, we move to the next case of homological mirror symmetry: elliptic curves.

For the purposes of mirror symmetry, an elliptic curve is a smooth Calabi-Yau variety of dimension one. We shall 
stick to $\C$ as a ground field. A smooth variety $X$ is Calabi-Yau if its canonical bundle is trivial. 
For most Calabi-Yau varieties, mirror symmetry is a little more symmetric. Unlike the previous section 
where we had four different definitions of categories of branes, for Calabi-Yau's we only have two. 
Given a Calabi-Yau $X$, the category of B-branes on $X$ is $D^b(\Coh(X))$ and the category of A-branes 
is $D^{\pi}(\Fuk(X))$. Homological mirror symmetry should simply exchange A-branes and B-branes of each ellipitic curve.

Testing this conjecture on the case of elliptic curves allows us to exploit some additional structure. 
We can uniformize elliptic curves. Namely, we can present any elliptic curve as $\C/(\Z+\tau\Z)$ for 
some $\tau$ in the upper half of the complex plane. This determines its complex structure, as we simply 
take the complex structure provided by $\C$. This also determines the symplectic structure, since 
$dz\wedge d\bar{z}$ is invariant under translations. We can of course scale the volume of our elliptic 
curve without changing the complex structure. This makes sense for any positive real number, but we will 
want to complexify and allow scalings by complex numbers with positive imaginary components. The resulting 
complex-valued symplectic form is usually called the complexified K\"ahler form. As such, we can describe 
our elliptic curve $E$ with complexified K\"ahler form by two numbers $\tau$ and $\rho$ which lie in the 
upper half plane. Let us incorporate this observation by denoting our curve by $E^{\rho}_{\tau}$. Then, 
mirror symmetry is simply the exchange $E^{\rho}_{\tau} \leftrightarrow E^{\tau}_{\rho}$.

Below, we generally follow the notation from \cite{PZ98}.

\subsection{A-branes on an elliptic curve}

Another useful property of elliptic curves is that mean curvature flow behaves well on them. One 
can also note that mean curvature flow is a Hamiltonian deformation \cite{TY02}. Thus, we can flow 
our (homologically non-trivial) Lagrangian (curve) into a geodesic and be none the worse off in terms 
of the Fukaya category. The result is that we only need to consider special Lagrangians when working 
with the Fukaya category. Of course, the standard means of computing $CF(L,L)$ is by using a slight 
Hamiltonian deformation so we need either to modify our definitions or forget about these types of 
computations. While one could attempt generalisations of Fukaya-type compositions using singular 
homology and incidence requirements \cite{FOOO}, such a discussion would certainly take us beyond 
the scope of our introductory paper. We shall therefore restrict our attention to a subclass of 
morphisms when defining compositions.

We reduce to the following. The objects of our category will be special Lagrangian submanifolds of $E^{\rho}_1$. 
Our special Lagrangians are simply lines with rational (or infinite) slope. We will of course need 
to choose a grading and a spin structure as before. The choice of spin structure we shall suppress. 
Let us recall how a grading of a Lagrangian $L$ is chosen. We need to lift the phase map $L \ra S^1$. 
Our holomorphic volume form is simply $dz$. Restricting that volume form to the line $mx=ny$, we get 
$\frac{(n+\imath m)}{\sqrt{m^2+n^2}}dl$ where $l$ is some coordinate and we have 
$l \ra (1/\sqrt{m^2+n^2})(nl,ml)$. With the parametrization, the induced volume form on our line is 
simply $dl$. The phase map is constant and equal to $e^{2\pi\imath \a}$ for some choice of $\a$ so that 
$e^{\pi \imath \a}$ is equal to the argument of $n+\imath m$. Thus, to choose a grading for our 
Lagrangian $L$ we simply need to make a choice of $\a$. There are $\Z$-fold options and each choice 
gives a different object in our category of A-branes.

Now, the grading (which we take from \cite{PZ98}) on an intersection point $p \in L_i \cap L_j$ is 
then $\mu(p) = -[\a_j-\a_i]$ where $[r]$ is the greatest integer less than or equal to $r$. If 
we exclude vertical lines, we can choose a grading that satisfies 
$\mu(p) = 0$ if $s_i < s_j$ and $\mu(p) = 1$ if $s_i > s_j$ for all lines simultaneously.

We also include flat line bundles $\mathcal{E}$ on our special Lagrangians with
monodromy of unit modulus (or flat $U(1)$-bundles). Then an A-brane on our elliptic curve is a 
triple $(L,\a,\mathcal{E})$ with $L$ a special Lagrangian, $\a$ a lift of the argument of the slope, 
and $\mathcal{E}$ a flat line bundle on $L$. We shall denote the collection of this data as $\mathcal{L}$.

Let us set $\Hom(\mathcal{L}_i,\mathcal{L}_j) = CF(L_i,L_j; \Hom(\mathcal{E}_i,
\mathcal{E}_j))$. 

As discussed above, we restrict the class of allowed compositions. In our
$\Ainf$-category, we are only allowed to form compositions 
$$m_n:\Hom(\mathcal{L}_1,\mathcal{L}_2)\otimes\cdots\otimes\Hom(\mathcal{L}_{n-1},L_n) \ra 
\Hom(\mathcal{L}_1,\mathcal{L}_n)$$ where all Lagrangians involved are
distinct. This of course forces us to abandon identity morphisms.

We have the compositions.
\begin{displaymath}
 m_n:\Hom(\mathcal{L}_1,\mathcal{L}_2)\otimes\cdots\otimes\Hom(\mathcal{L}_{n-1},L_n)
\ra \Hom(\mathcal{L}_1,\mathcal{L}_n)
\end{displaymath}
\begin{displaymath}
 m_n(p_1,\cdots,p_n) = \sum_{p_0 \in L_1 \cap L_n} \sum_{\phi \in \mathcal{M}(p_0,\cdots, p_n)} \pm \exp(2\pi\imath\int_D \phi^*\omega) P\exp(\phi^*\beta)
\end{displaymath}
when $\deg(p_0) = \sum_{i=1}^n \deg(p_i) + 2-n$, $\mathcal{M}(p_0,\cdots,p_n)$ is the
space of holomorphic maps from an $(n+1)$-pointed disc (modulo automorphisms when $n=1$) with 
$\phi(\xi_i) = p_i$ and $\phi(S_i) \subset L_i$, and $P\exp(\phi^*\beta)$ is the path ordered exponential 
\begin{displaymath}
 P\exp(\phi^*\beta) = P\exp(\int_{S_n} \beta_n) t_n P\exp(\int_{S_{n-1}} \beta_n)
t_{n-1} \cdots t_1 P\exp(\int_{S_1} \beta_n)
\end{displaymath}
Here $\beta$ is a connection form for the flat bundle $\mathcal{E}$.

Consider two A-branes $\mathcal{L}$ and $\mathcal{L}'$. Then as a set the morphisms
between $\mathcal{L}$ and $\mathcal{L}'$ are the same --- simply $L \cap L'$. Assume we
have equipped these A-branes with a grading with phase within $(-1/2,1/2)$. Then,
all the morphisms between $\mathcal{L}$ and $\mathcal{L}'$ are either of degree zero
or of degree one. Let us assume they are of degree zero. Assume that $L$ has slope
$m$ and $L'$ has slope $n$. Then there are $n-m$ intersection points on the torus
given by
\begin{displaymath}
 \left(\frac{k}{n-m},\frac{mk}{n-m}\right)
\end{displaymath}
for $k \in \Z/(n-m)\Z$.
Assume now that the $x$-axis intercept of $L$ is $(\a_1,0)$ and the $x$-axis intercept for
$L'$ is $\a_2$ for $0\leq \a_i \leq 1$. Then the equation of $L$ becomes
$(\a_1+t,(n-1)\a_1+nt)$ and the equation of $L'$ is $(\a_2+t,(m-1)\a_2+mt)$. Thus
the intersection points are
\begin{displaymath}
 \left(\frac{k+\a_2-\a_1}{n-m},\frac{mk+n\a_2-n\a_1}{n-m}\right)
\end{displaymath}
for $k \in \Z/(n-m)\Z$.

There is a natural non-degenerate pairing of degree one
\begin{displaymath}
 \Hom(\mathcal{L},\mathcal{L'}) \times \Hom(\mathcal{L'},\mathcal{L}) \ra \C
\end{displaymath}
where we declare $p \in \Hom(\mathcal{L},\mathcal{L'})$ dual to $p \in
\Hom(\mathcal{L'},\mathcal{L})$. We then take the trace on the bundle components.
Let us denote this by $(\cdot,\cdot)$. Then
\begin{displaymath}
 (m_2(a,b),c) = (a,m_2(b,c))
\end{displaymath}
If we set $\mathcal{L} = \mathcal{L}'$, then the pairing $(\cdot,\cdot)$ is the Floer
manifestation of Poincar\'e duality. The cyclic symmetry and nondegeneracy of the
pairing allows us to reconstruct two-fold compositions involving degree one morphisms 
from the two-fold compositions only involving degree zero morphisms. Thus, at the level of 
graded categories, if we formally add in identity morphisms to $H^0(\Fuk(E^{\tau}))$ then we 
can determine the proper two-fold composition rules for the degree one piece of 
$\Hom(\mathcal{L},\mathcal{L})$.

\subsection{B-branes on an elliptic curve}
We now move onto the other side. If we quotient $\C$ out by $z \mapsto z+1$ first, our elliptic curve $E_{\tau}$ can be viewed as $\C^{\times}/\Z$ where the action by $\Z$ is given by multiplication by $e^{2\pi\imath\tau}$. We can pull back any line bundle from $E_{\tau}$ to $\C^{\times}$. Any line
bundle on $\C^{\times}$ is necessarily trivial (and consequently so is the pullback of
any vector bundle). Since $E_{\tau}$ is projective and smooth, the collection of all
line bundles on $E_{\tau}$ generates $D^b(\Coh(E_{\tau}))$. All line bundles, since they are
trivial when pulled back to $\C^{\times}$, can be described as the quotient
$(\C^{\times} \times \C)/\Z$ where $\Z$ acts by $(u,v) \mapsto (qu,\phi(u)v)$ for
some holomorphic map $\phi: \C^{\times} \ra \C^{\times}$. The map $\phi$ is the choice of
holomorphic trivialisation over $\C^{\times}$ of the pulled-back bundle. Let $L(\phi)$
denote the line bundle determined in this way.

With the similarity of the Fukaya-compositions to theta functions kept in mind, there
are obvious choices for $\phi$. Let $\phi_0(u) = \exp(-\pi\imath\tau)u^{-1}$. The
sections of $L(\phi_0)$ are the classical theta functions. A general theta function is
given as follows
\begin{displaymath}
 \theta[c',c''](\tau,z) = \sum_{m\in\Z}\exp(2\pi\imath(\tau(m+c')^2/2+
(m+c')(z+c'')))
\end{displaymath}
The classical theta function is simply $\theta[0,0](\tau,z)$. To check that this is actually
a section of $L$, we simply need to check that it behaves properly when we translate by
$1$ and $\tau$. $\theta[0,0](\tau,z)$ is clearly invariant under the shift $z \mapsto
z+1$. And
\begin{displaymath}
 \theta[0,0](\tau,z+\tau) = \sum_{m=-\infty}^{\infty} \exp(2\pi\imath((m^2/2)\tau + mz + m\tau))
\end{displaymath}
Letting $l = m+1$, we have $(m^2/2)\tau + mz + m\tau = (l^2/2)\tau+lz-(1/2)\tau-z$. Thus,
\begin{displaymath}
 \theta[0,0](\tau,z+\tau) = e^{-\pi \imath \tau}e^{-2\pi \imath z}\theta[0,0](\tau,z),
\end{displaymath}
so it is indeed a section of $L$. One can similarly check that $\theta[a/n,0](n\tau,nz)$
for $a \in \Z/n\Z$ are sections of $L^n$. By pulling back $L^n$ via translation, we get
sections for all possible line bundles.

The degree of $L(\phi_0) = L$ is $1$. Any line bundle on an elliptic curve is the tensor
product of a degree zero line bundle and $L^n$. It is well known that $\Pic^0(E_{\tau})$
is isomorphic to $E_{\tau}$. Thus, any line bundle is simply $t_x^*L \otimes L^{n-1}$ where
$x \in E$ and $n$ is the degree. If the degree of a line bundle is $< 0$,
the corresponding divisor is not equivalent
to an effective divisor. Hence, we can have no global sections.
For $n>0$, the sections of $L^n$ are also given by theta functions:
$\theta[a/n,0](n\tau,nz)$ for $a\in \Z/n\Z$. The sections of $t_{\a+\imath \b}^*L^n$ are then
\begin{displaymath}
 t_{\a + \imath \b}^* \theta[a/n,0](n\tau,nz) = \theta[a/n,0](n\tau,n(z+\a+\imath\b))
\end{displaymath}
Note that
\begin{displaymath}
 \Hom(t_x^*L \otimes L^{n-1},t_y^*L \otimes L^{m-1}) \cong \Hom(t_x^*L^{-1}\otimes t_y^*L \otimes L^{m-n})
\end{displaymath}
Similarly, $t_x^*L^{-1} \otimes t_y^*L^{-1}$ can have no global sections if $x \not =
y$. Thus, there are no morphisms from $t_x^*L \otimes L^{n-1}$ to $t_y^*L \otimes L^{m-1}$ unless 
$m > n$ or $m = n$ and $x = y$. If we let $x = \a_1 + \imath \b_1$ and $y = \a_2 + \imath b_2$, then 
we can also rewrite $t_x^*L^{-1} \otimes t_y^*L^{-1}$ as $(t_{\a_{21}+\imath\b_{21}}^*L)^{m-n}$ where
\begin{displaymath}
 \a_{21} = \frac{\a_2 - \a_1}{m-n}, \b_{21} = \frac{\b_2-\b_1}{m-n}
\end{displaymath}
The sections of $(t_{\a_{21}+\imath\b_{21}}^*L)^{m-n}$ are given by
\begin{displaymath}
 \theta[a/(m-n),0]((m-n)\tau,(m-n)(z+\a_{21}+\imath\b_{21}))
\end{displaymath}

Since $E_{\tau}$ is a Calabi-Yau variety, the canonical bundle $\omega_{E_{\tau}}$ is
trivial. For a general smooth projective $X$, we have $\Ext^*(F,G) \cong
\Ext^{n-*}(G,F\otimes \omega_X)^{\vee}$
\cite{Hart77}. The pairing is given by composing $a \in \Ext^*(F,G)$ and $b \in
\Ext^{n-*}(G,F\otimes \omega_X)$, taking the trace on $F$ giving an element of
$H^n(\omega_X) \cong \C$. The functor $S: D^b(\Coh(X)) \ra D^b(\Coh(X))$ given by 
$F \mapsto F \otimes \omega_X[-\dim X]$ is called the Serre functor \cite{BK89}. 
For a Calabi-Yau, we get something rather special as $S \cong [-\dim X]$. We also have 
\begin{displaymath}
 (a\cdot b,c) = (-1)^{\tilde{a}(\tilde{b}+\tilde{c})} (b \cdot c,a)
\end{displaymath}
Thus, as in the case of A-branes, one can then use Serre duality and the degree zero
morphisms to compute the Ext-groups of these line bundles and reconstruct
$D^b(\Coh(E_{\tau}))$ as a graded category.

In \cite{PZ98}, Polishchuk and Zaslow treat all possible vector bundles on $E_{\tau}$ 
and torsion sheaves.

\subsection{An outline of the equivalence}

The first step in checking homological mirror symmetry for elliptic curves is
to give a functor from the category of line bundles (with only honest morphisms)
to the sub-category of the Fukaya category given by lines of non-vertical slope
with flat $U(1)$-bundles on them. From the discussion above, one can see that
the correspondence between theta functions and intersection points is pretty
strong, so we make the obvious definition. Let
$$\Phi: \mathcal{L}(E_{\tau}) \ra H^0(\Fuk(E^{\tau}))$$
$$\Phi(t_{\a \tau + \beta}^*L\otimes L^{n-1}) = (\Lambda,A)$$
where $\Lambda = (\a + t,(n-1)\a + nt)$ and $A = (-2 \pi \imath \beta)dx$.
Recall that morphisms between $t_{\a_1 \tau + \beta_1}^*L\otimes L^{n_1-1}$
and $t_{\a_2 \tau + \beta_2}^*L\otimes L^{n_2-1}$ are the same thing as sections
of $t_{\a_{21}\tau+\b_{21}}^*L^{n_2-n_1}$ where
$$\a_{21} = \frac{\a_2-\a_1}{n_2-n_1}, \b_{21} = \frac{\b_2-\b_1}{n_2-n_1}.$$
$t_{\a_{21}\tau+\b_{21}}^*L^{n_2-n_1}$ has sections given by
$$f_k\theta[k/(n_2-n_1),0]((n_2-n_1)\tau,(z+\a_{21}\tau + \b_{21}))$$
for $k \in \Z/(n_2-n_1)\Z$. Correspondingly the intersection points of
$$\Phi(t_{\a_1 \tau + \beta_1}^*L\otimes L^{n_1-1})\,\,\mbox{and}\,\, \Phi(t_{\a_2 \tau +
\beta_2}^*L\otimes L^{n_1-1})$$ are given by 
$$e_k = \left(\frac{k+\a_2-\a_1}{n_2-n_1},\frac{n_1k+n_1\a_2-n_1\a_1}{n_2-n_1}
\right)$$
for $k\in\Z/(n_2-n_1)\Z$.
We then set $\Phi(f_k) = \exp(-\pi \imath \tau \a_{21}^2(n_2-n_1)) \cdot e_k$.

\begin{thm}
 \cite{PZ98} $\Phi$ is full and faithful. It gives an equivalence with the full subcategory
of $H^0(\Fuk(E^{\tau}))$ consisting of special Lagrangians with non-vertical slope
and the chosen special gradings paired with flat $U(1)$-bundles. 
Moreover, $\Phi$ naturally extends to an equivalence
of the abelian category of coherent sheaves on $E_{\tau}$ and 
the subcategory $H^0(\Fuk(E^{\tau}))$ consisting of all Lagrangians with flat
$U(N)$-bundles and the chosen special gradings.
\end{thm}

In lieu of reviewing the proof, we shall treat an illuminating example. The example
is taken directly from \cite{PZ98}, as the author feels he cannot improve upon the
choice.

Under $\Phi$ we take the structure sheaf of $E_{\tau}$ to the $x$-axis with zero connection.
$L^n$ is mapped to a line passing through the origin with slope $n$ and no connection. Let us
consider the composition $\Hom(\mathcal{O},L) \times \Hom(L,L^2) \ra \Hom(\mathcal{O},L^2)$.
The vector space $\Hom(\mathcal{O},L)$ is one-dimensional and spanned by the classical
theta function $\theta[0,0](\tau,z)$. As $\Hom(L,L^2) \cong \Hom(\mathcal{O},L)$, this
space is also spanned by the classical theta function.

On the other side, we have $\Lambda_1 \cap \Lambda_2 = \{e_1\}, \Lambda_2 \cap \Lambda_3
= \{e_1\}$, and $\Lambda_1 \cap \Lambda_3 = \{e_2,e_3\}$ where $e_1$ and $e_2$ are
the origin and $e_3 = (1/2,0)$. Under $\Phi$, we identify $e_1$ with $\theta[0,0](\tau,z)$
and $e_2$ with $\theta[0,0](2\tau,2z)$ and $e_3$ with $\theta[1/2,0](2\tau,2z)$.

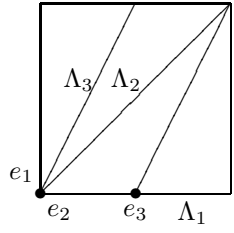
\begin{displaymath}
 \xy
 0;/r.30pc/:
 (-10,-10)*{}; (10,-10)*{} **\dir{-};
 (-10,-10)*{}; (-10,10)*{} **\dir{-};
 (10,-10)*{}; (10,10)*{} **\dir{-};
 (-10,10)*{}; (10,10)*{} **\dir{-};
 (-10,-10)*{}; (10,10)*{} **\dir{-};
 (-10,-10)*{}; (0,10)*{} **\dir{-};
 (0,-10)*{}; (10,10)*{} **\dir{-};
 (6,-12)*+{\Lambda_1};
 (-1,2)*+{\Lambda_2};
 (-6,2)*+{\Lambda_3};
 (-12,-8)*+{e_1};
 (-8,-12)*+{e_2};
 (-10,-10)*+{\bullet};
 (0,-10)*+{\bullet};
 (0,-12)*+{e_3};
 \endxy
\end{displaymath}

To compute the composition in the Fukaya category, we lift our diagram to $\R^2$ in all
possible ways and count the number of triangles with appropriate weights determined by
the K\"ahler form. These will give the structure coefficients.

Of course, triangles which are translates of another triangle are equivalent. Therefore,
we need to pin down one of the vertices that we will take to be the origin, i.e. $e_1$.
We need to compute the coefficients of $e_2$ and $e_3$. The coefficient of $e_2$ comes
from counting triangles with vertices on integer translates of the origin. Since we have 
pinned down the first vertex we are looking at triangles in $\R^2$ given by
$(0,0),(n,0),(2n,2n)$. Each such triangle has area $\tau n^2$. Summing over them we get
$$\sum_{n=-\infty}^{\infty} \exp(2\pi\imath\tau n^2) = \theta[0,0](2\tau,0)$$
The coefficient of $e_3$ comes from triangles with vertices at $$(0,0),(0,n+1/2),(2n+1,
2n+1).$$ We get
$$\sum_{n=-\infty}^{\infty} \exp(2\pi\imath\tau(n+1/2)^2) = \theta[1/2,0](2\tau,0).$$
Thus, we have
$$m_2(e_1,e_1) = \theta[0,0](2\tau,0)e_2 + \theta[1/2,0](2\tau,0)e_3.$$

Let us now investigate the other side. We are looking at the square of the classical
theta function. Since $\Phi$ is a functor, we need to show that
$$\theta(\tau,z)^2 = \theta[0,0](2\tau,0)\theta[0,0](2\tau,2z) + \theta[1/2,0](2\tau,0)
\theta[1/2,0](2\tau,2z)$$
This is simply an application of the addition formula for theta functions
\cite{Mum84}. This is a basic example but it nicely illustrates the correspondence. 

Using the equivalence between $\Coh(E_{\tau})$ and the subcategory of $H^0(\Fuk(E^{\tau}))$ with objects whose grading is in the interval $(-1/2,1/2)$, 
one can then extend the functor to $\Phi: D^b(\Coh(E_{\tau})) \ra H^0(\Fuk(E^{\tau}))$ uniquely by
requiring that it be compatible with Serre and Poincar\'e duality. Thus, we get

\begin{thm}
 $D^b(\Coh(E_{\tau}))$ and $H^0(\Fuk(E^{\tau}))$ are equivalent as graded categories.
\end{thm}

The next step is to determine if the triangulated structures coincide. The triangulated
structure descends from a cone construction at the chain level. Thus, the natural
way to proceed is to show that as $\Ainf$-categories $D^b(E_{\tau})$ and $H^0(\Fuk(E^{\tau}))$
are quasi-isomorphic. But recall that the definition of the Fukaya category is
missing $\Hom(\mathcal{L},\mathcal{L})$. We must restrict our class of composable morphisms
to Hom-spaces between different objects. With this limitation in mind, one can then
attack the problem. 

In \cite{Pol04}, the possible minimal $\Ainf$-structures on $D^b(\Coh(E_{\tau}))$ compatible
with Serre duality (i.e. cyclic) and the tensor structure and with a suitably large
class of allowed compositions are classified up to homotopy. The problem then
reduces to computing a single triple product as the test of homotopy equivalence. The
triple products on $D^b(\Coh(E_{\tau}))$ and $H^0(\Fuk(E^{\tau}))$ coincide. Thus, we get
the following result.

\begin{thm}
 Consider $D^b(\Coh(E_{\tau}))$ with the admissible class of morphisms determined by the
equivalence with $H^0(\Fuk(E^{\tau}))$. The $\Ainf$-structures induced on this category
by the dg-category of line bundles and the Fukaya category coincide.
\end{thm}

In fact, the result is a little stronger. If one can extend the definition of the Fukaya
category to include endomorphism sets while preserving the $\Ainf$-structure already
present on composable morphisms, then one knows that the resulting $\Ainf$-category will
be homotopy equivalent to $D^b(E_{\tau})$ with the minimal $\Ainf$-structure determined
by the dg-category of line bundles. Of course, one can extend the $\Ainf$-structure in
such a way since the $\Ainf$-structure inherited from the dg-category of line bundles
allows all possible compositions. So such an extension certainly exists; it is just not
of a symplectic nature yet. Modulo this small detail we have homological mirror symmetry
for elliptic curves.

\section{Further results}
\label{sec:furtherresults}

In this section, we give a rapid review of the established cases of homological mirror 
symmetry beyond dimension one. By no means is this meant to be a comprehensive treatment. 
The interested reader is encouraged to consult the references.

\subsection{A-branes on (near) Fano manifolds versus B-branes on LG-models}

\subsubsection{Fano surfaces}

In the previous examples, we were fortunate to be able to compute compositions in
Fukaya categories. The Landau-Ginzburg mirror to $\P_k^1$ was an extremely simple situation
since the fibers of the potential were zero-dimensional manifolds as were the Lagrangians
involved. For $\P_{\C}^1$, any choice of almost complex structure was automatically integrable
and regular. Thus, we could apply the Riemann mapping theorem and reduce our counts to
purely topological considerations.

If we consider A-branes on the mirror of $\P^2$, the fibers are no longer zero-dimensional 
but are now non-compact Riemann surfaces. Thus, the counting can again be reduced to topological 
considerations. This at least gives one the hope of applying previous methods to this problem. 
Moreover, if we replace $\P^2$ by any other two-dimensional smooth variety (or stack), the 
Landau-Ginzburg mirror (if it exists) should be approachable via the same methods.

In fact, this is done in a pair of papers \cite{AKO04,AKO06}, see also \cite{Ued06}. 
In \cite{AKO04}, Auroux, Katzarkov, and Orlov tackle weighted projective surfaces as stacks and demonstrate 
the following.

\begin{thm}
 (Auroux-Katzarkov-Orlov) Let $a,b,c$ be mutually prime positive integers and consider 
the stack $\C\P(a,b,c)$. Then $D^b(\Coh(\C\P(a,b,c))$ is triangle equivalent to  
$D^{\pi}(FS(W))$ where $W: \lbrace x^ay^bz^c = 1 | x,y,z \in \C \rbrace \ra \C$ is given 
by $W(x,y,z) = x+y+z$.
\end{thm}

Moreover, they also address how homological mirror symmetry behaves with respect to 
non-commutative deformations of the homogeneous coordinate ring of the weighted projective plane. 
In addition, using the results for weighted projective planes, they prove a similar result for 
Hirzebruch surfaces.

In \cite{AKO06}, Auroux, Katzarkov, and Orlov continue the study of this side of mirror symmetry for surfaces. 
They investigate the effect of blowing up $\C\P^2$ at $\leq 9$ points. To properly address this, 
they must compactify the mirror to $\C\P^2$ and deform it.

\begin{thm}
 (Auroux-Katzarkov-Orlov) Let $K$ be a set of $k$ points in $\C\P^2$. Consider the del Pezzo 
surface $X_K$ given by blowing up $\C\P^2$ at those points. Then $D^b(\Coh(X_K))$ is triangle 
equivalent to $D^{\pi}(FS(W_K))$ where $W_K : M_K \ra \C$ is an elliptic fibration obtained by 
deforming a compactified mirror of $\C\P^2$.
\end{thm}

For more details, see \cite{AKO04,AKO06}. In \cite{Ued06}, Ueda studies toric blow-ups of $\C\P^2$.

\subsubsection{Projective, smooth toric varieties}

In \cite{Abou07}, Abouzaid takes a slightly different approach to the study of homological 
mirror symmetry for smooth projective toric varieties. Motivated by the results on Fano surfaces 
and the underlying physics, one can speculate that the mirror of a such a toric variety $X_{\Delta}$ 
determined by the polytope $\Delta$ is the Landau-Ginzburg model 
$W_{\Delta}: (\C^{\times})^{\dim X_{\Delta}} \ra \C$ where $W$ is the Laurent polynomial whose Newton 
polytope is $\Delta$. In fact, symplectically, one has a lot of freedom in choosing the coefficients 
of $W_{\Delta}$. This allows one to take a limit that corresponds to a tropical degeneration. In this 
limit, one can identify A and B branes on either side with certain Morse chain complexes on the 
moment polytope of $X_{\Delta}$. Morally, one reduces mirror symmetry in this setting to data on the 
shared moment polytope. Of course, a bit of work goes into relating these ``tropical Lagrangians'' back 
to the usual Lagrangians. In the end, one gets the following result. Let $\Fuk(W_{\Delta})$ denote 
the $\Ainf$-(pre)-category of admissible Lagrangians in $(\C^{\times})^{\dim X_{\Delta}}$. Essentially, 
these are Lagrangians with boundary on a regular fiber over a point $p$ of $W_{\Delta}$ which project 
under $W_{\Delta}$, near $p$, to curves emanating from $p$.

\begin{thm}
 (Abouzaid) There is a full and faithful triangle functor 
$$i: D^b(\Coh(X_{\Delta})) \ra D^{\pi}(\Fuk(W_{\Delta})).$$ Moreover, 
if $X_{\Delta}$ is ample, the functor is essentially surjective.
\end{thm}

For more details, see \cite{Abou06}.

\subsection{A-branes on Fano manifolds versus B-branes on LG-models}

Although there are few results on the Fukaya category of a toric variety, similar to subsection 
\ref{subsec:A-branesonP1}, one can define a symplectic category that works. In \cite{CO06}, 
Cho and Oh define a version of Floer cohomology called adapted Floer cohomology which uses the 
standard complex structure on the toric variety instead of a generic almost complex structure. If 
the toric variety is convex, then this reduces to the usual Floer cohomology. They also restrict 
themselves to the Lagrangian torus fibers (with flat line bundles) of the moment map fibration.

As in subsection \ref{subsec:A-branesonP1}, there are only a finite number of torus fibers that 
have non-vanishing Floer cohomology. These are called balanced torus fibers. Cho and Oh show 
that balanced torus fibers are in one-to-one correspondence with the critical points of the mirror 
LG superpotential. There is no Floer cohomology between two distinct balanced torus fibers. And, 
given a balanced torus fiber $L$, the $\Z/2\Z$-graded Floer cohomology is $HF(L,L)$ is a Clifford 
algebra on $\dim L$ generators with non-degenerate pairing given by the holomorphic Hessian of $W$ 
at the critical point $p$ corresponding to $L$. As in subsection \ref{subsec:A-branesonP1}, the 
intrinsic formality of the Clifford algebra uniquely determines any underlying $\Ainf$-structure on $HF(L,L)$.

For B-branes on the mirror, a simple computation shows that the skyscraper sheaves at the singular 
points are mutually orthogonal, and that the endomorphism algebra of any of these objects is a 
Clifford algebra on $\dim L$ generators with the bilinear form given by the holomorphic Hessian at 
the corresponding point. As before, to match up completely, we need to change variables from 
$\C^{\times}$ to $\C$ using the exponential map. Nevertheless, we get the following.

\begin{thm}
 For a convex smooth projective toric variety $X_{\Delta}$ with mirror $W_{\Delta}$, the smallest 
triangulated thick subcategory contaning the skyscraper sheaves at the singular points 
of $W_{\Delta}$ in $D_{sing}(W_{\Delta})$ 
is triangle equivalent to the smallest triangulated thick subcategory containing 
the balanced torus fibers in $D^{\pi}(\Fuk(X_{\Delta}))$.
\end{thm}

\subsection{A-branes versus B-branes on Calabi-Yau manifolds}

\subsubsection{Abelian varieties}

Following the results of \cite{PZ98}, the natural inclination was to consider higher dimensional 
abelian varieties $A$ using similar methods. And, indeed, in \cite{Fuk02}, this is what is done. 
Fukaya is able to demonstrate that the corresponding Hom-spaces match up and, in the case 
where all three Lagrangians are mutually transverse, the $m_3$-compositions match also.

Another approach was undertaken in \cite{KS01}. Kontsevich and Soibelman use non-Archimedean analysis to skirt 
convergence issues and, similar to \cite{Abou07}, use a degeneration to the base of a torus fibration 
to show that $D^b(\Coh(A))$ embeds as a triangulated subcategory of 
$D^{\pi}(\Fuk(A^{\vee}))$ where $A^{\vee}$ is the mirror abelian variety to $A$. Specifically, 
it is equivalent to the category of Lagrangians transverse to the base.

\subsubsection{Quartic surface}

In \cite{Sei03}, Seidel implements his plan from \cite{Sei02} which is to use deformation 
theory and results about directed categories of vanishing cycles to deduce homological mirror symmetry 
for a Calabi-Yau variety. The idea is as follows. The computations on the algebro-geometric side are 
not so bad. On the symplecto-geometric side, all the relevant Lagrangians sit in an affine patch of 
the Calabi-Yau, in this case a quartic surface in $\P^3$, as vanishing cycles of a Picard-Lefschetz 
fibration. Here they can be viewed as matching paths corresponding to lower-dimensional Picard-Lefschetz 
fibrations. This viewpoint allows one to automate an inductive procedure that computes Floer homology 
(Fukaya categories) associated to these Lagrangians from the directed category of vanishing cycles 
associated to the Picard-Lefschetz fibration. Then, to pass from the Lagrangians in the affine patch 
to the Lagrangians in the total space one has to include pseudo-holomorphic discs that hit the divisor 
at infinity, i.e. the complement of the affine patch. Counting these discs can be viewed as a formal 
deformation of the $\Ainf$-structure associated endomorphism algebra of the vanishing cycles in the 
affine patch. Then, the endomorphism algebra of the vanishing cycles in the compact Calabi-Yau is obtained by
specialising this formal deformation to a non-zero value of the parameter. The key to making this result 
feasible is that the deformation spaces involved are small, i.e. one-dimensional, and carry a 
$k^{\times}$ action. Thus, there are really only two different deformations possible, the trivial one 
and the non-trivial one. The issue of triviality of the deformation can be reduced to a computation 
involving Hochschild cohomology of the underlying $\Ainf$-algebra which can be checked on one side via 
explicit computations and on the other by the non-degeneration of a certain spectral sequence. Of course, 
Seidel still needs to avoid the convergence problem, so he works with Novikov rings.

\begin{thm}
 (Seidel) Let $X$ be a smooth quartic surface in $\P_{\C}^3$ and let $\Lambda_{\Q}$ be the rational 
Novikov field over $\C$. Take the quartic surface in $\P^3_{\Lambda_{\Q}}$ defined by
\begin{displaymath}
 x^4+y^4+z^4+w^4+q(xyzw) = 0,
\end{displaymath}
quotient by the standard action of $(\Z/4\Z)^2$, take the minimal crepant resolution, and denote it 
by $Y$. Then there is an equivalence of triangulated categories 
$$D^{\pi}(\Fuk(X)) \cong \hat{\psi}^*D^b(\Coh(Y))$$ 
where $\psi$ is a continuous automorphism of $\Lambda_{\N}$.
\end{thm}

For more details, see \cite{Sei03}.

\section{How is homological mirror symmetry related to mirror symmetry?}
\label{sec:HMSMS}

Mirror symmetry is a wide and varied subject. Homological mirror symmetry is a single attempt to 
unify a majority of the various topics. As such, it must prove its utility in other mirror symmetric 
mathematical investigations. In this section, we present the relations between homological 
mirror symmetry and three other prominent mathematical aspects of mirror symmetry.

\subsection{SYZ}

Perhaps the most tantalising approach to mirror symmetry was outlined by Strominger, Yau, and 
Zaslow in \cite{SYZ96}. Here, mirror symmetry of two Calabi-Yau manifolds could be realised 
geometrically, as opposed to homologically. Namely, given a Calabi-Yau manifold $X$ over $\C$ 
there should exist a fibration by special Lagrangian tori. Dualising this torus fibration 
(and accounting for the singular fibers appropriately) should realise the mirror Calabi-Yau $Y$. 
Such an explicit geometric operation should provide an explicit functor between the categories of 
branes on each side of the mirror correspondence. Unfortunately, special Lagrangian torus fibrations 
are hard to come by on Calabi-Yau manifolds. Of course, many expect mirror symmetry to hold not 
necessarily throughout the whole moduli spaces of theories on each side, but in a small neighborhood 
of some special singularities on each side. Thus, one can try to implement SYZ in conjunction with 
a degeneration of the structures. Indeed, motivated by the SYZ conjecture,
M. Gross and B. Siebert have begun a program \cite{GS06, GS07a, GS07b} to
explain mirror symmetry as a Legendre transform of dual affine structures obtained from toric
degenerations. For an introduction see \cite{Gro08}.  We have already seen such logic employed in the previous sections 
in investigating homological mirror symmetry both in the case of Calabi-Yau varieties, specifically 
abelian varieties, and also in the case of Fano varieties. Indeed, perhaps the most spectacular 
application of SYZ-inspired ideas to homological mirror symmetry is found in \cite{Abou07}, where 
one can actually see how the branes on each side correspond by repackaging them as data on the base 
moment polytope. And SYZ informs our understanding of the other side of homological mirror symmetry 
for Fano varieties. It is under SYZ that the balanced torus fibers of \cite{CO06} correspond to 
the critical points of the mirror superpotential. This is explored in \cite{Hori02}.

\subsection{Mirror maps}

Another common area of investigation in mirror symmetry is the determination of the mirror map. 
The mirror map is an isomorphism between neighborhoods of two special singular points in moduli on 
each side of the mirror correspondence. Aside from its evident importance, mirror maps turn out to 
possess rich arithmetic properties. Mirror maps are interesting and elusive. If indeed homological 
mirror symmetry provides the fullest explanation for the mirror symmetry phenomenon, then we 
should somehow be able to extract the mirror map. Indeed, in some cases, one can do exactly that.

In \cite{AZ06}, Aldi and Zaslow use the following idea of Seidel. A projective variety $X$ is determined 
by its homogeneous coordinate ring. The data of the homogeneous coordinate ring sits inside 
$D^b(\Coh(X))$ and can be extracted assuming we know the very ample line bundle $L$ corresponding 
to the embedding and the auto-equivalence corresponding to $-\otimes L$. Then, the homogeneous coordinate 
ring is $\oplus_{n > 0} \Hom_{D^b(\Coh(X))}(L,L^{\otimes n})$. Given a triangulated category, an object, 
and an autoequivalence we can always form such a ring but there is no reason to expect commutativity. 
But, if we know that homological mirror symmetry holds, then we know what object and functor to look at --- 
the mirror to $L$ and the mirror to $-\otimes L$. In general, it is expected that the mirror autoequivalence 
to tensoring by a line bundle $L$ is given by the Dehn twist $\tau_S$ about a vanishing cycle $S$ mirror to 
$L$. Assuming homological mirror symmetry, we can extract the homogeneous coordinate ring of the mirror 
as $\oplus_{n > 0} \Hom_{D^{\pi}(\Fuk(Y))}(S,\tau^n(S))$. The structure of this ring depends only on the 
symplectic structure of $Y$. Thus, as we vary the symplectic structure of $Y$, we should see the 
variation of the complex structure of the mirror $X$ manifested in the change in the homogeneous coordinate 
ring of $X$. In this way, Aldi and Zaslow of \cite{AZ06} extract the mirror map in a number of examples 
and check that the logic is correct.

The upshot is that if one can check homological mirror symmetry in a given case by exhibiting an 
explicit equivalence, one should be able to compute the corresponding mirror map and extract the 
wealth of arithmetic information that is expected to lie within.

\subsection{Instanton numbers}

Mirror symmetry broke onto the mathematical scene by properly counting curves \cite{COGP91} on Calabi-Yau 
threefolds, and curve, or instanton, counts have been a central feature of mirror symmetry since. 
Indeed, mirror symmetry was originally viewed as the study of these counts and their relation 
to solutions of certain differential equations. If homological mirror symmetry provides the fullest 
explanation of the mirror symmetry phenomenon, then one should be able to extract these instanton counts 
solely from the categorical data. Physically, the instanton counts come from the closed string sector of 
our topological field theory, whereas what has been discussed earlier concerns the open string sector. 
The objects in each of our categories are appropriate boundary conditions for these open strings in 
different theories. One needs a mathematical construction that relates the two. In \cite{Cos05,Cos07}, 
Costello considers the notion of a open-closed topological conformal field theory inside which sit 
open and closed TCFTs. He shows that an open TCFT is (morally) the same as an $\Ainf$-category with a 
cyclically-invariant inner pairing and demonstrates a means to construct, from an open TCFT, an open-closed 
TCFT whose space of closed states is the Hochschild chain complex $C_*(A,A)$ of the $\Ainf$-category 
associated to the initial open TCFT. He also demonstrates how to build a Gromov-Witten potential from 
a TCFT and thus gives instanton counts associated to an $\Ainf$-category with a cyclically-invariant 
inner pairing. In the case of the Fukaya category, the passage to Hochschild homology has a more natural 
explanation. For an $\Ainf$-category with a cyclically invariant inner product, one has a quasi-isomorphism 
between $C_*(A,A)$ and $C^*(A,A)$, the Hochschild chain complex. The Hochschild cohomology of an 
$\Ainf$-category $A$ is the same as the space of endomorphisms of $A$ in the derived category of 
$A$-bimodules; here a bimodule over a category is simply an autofunctor of the category. Using the tensor 
product, we see that $A$ corresponds to the identity functor. For the Fukaya category of a symplectic 
manifold $(M,\omega)$, we can make autofunctors from Lagrangians in $(M \times M, -\omega \oplus \omega)$. 
Morphisms between these functors correspond to Floer cohomology between the associated Lagrangians. 
The identity corresponds simply to the diagonal and it is now well known that $HF(L,L)$ is isomorphic 
to the quantum cohomology of $L$. For more on Lagrangian correspondences, see \cite{WW07}.

\bigskip

\bibliographystyle{nhamsplain}
\bibliography{refup}

\end{document}